\documentclass[11pt, a4paper]{article}
\usepackage[cp1251]{inputenc}
\usepackage{amsmath} \usepackage{euscript} \usepackage{marvosym}
\usepackage[dvipsnames]{xcolor}
\usepackage{mymatrx6}%\let\mymatrixII\mymatrix
\usepackage{graphicx}
\usepackage{longtable}
\usepackage{mathrsfs}
\usepackage{amssymb, amscd}
\usepackage{comment}
\def\mymatrix{\MyMatrixwithdelims..}
\Linestrue\Autonumfalse
\usepackage{mathrsfs}
\begin{document}

\newcounter{bnomer} \newcounter{snomer}
\newcounter{bsnomer}
\setcounter{bnomer}{0}
\renewcommand{\thesnomer}{\thebnomer.\arabic{snomer}}
\renewcommand{\thebsnomer}{\thebnomer.\arabic{bsnomer}}
\renewcommand{\refname}{\begin{center}\large{\textbf{References}}\end{center}}

\setcounter{MaxMatrixCols}{14}

\newcommand\restr[2]{{% we make the whole thing an ordinary symbol
  \left.\kern-\nulldelimiterspace % automatically resize the bar with \right
  #1 % the function
  %\vphantom{\big|} % pretend it's a little taller at normal size
  \right|_{#2} % this is the delimiter
}}

\newcommand{\sect}[1]{%
\setcounter{snomer}{0}\setcounter{bsnomer}{0}
\refstepcounter{bnomer}
\par\bigskip\begin{center}\large{\textbf{\arabic{bnomer}. {#1}}}\end{center}}
\newcommand{\sst}[1]{%
\refstepcounter{bsnomer}
\par\bigskip\textbf{\arabic{bnomer}.\arabic{bsnomer}. {#1}}\par}
\newcommand{\defi}[1]{%
\refstepcounter{snomer}
\par\medskip\textbf{Definition \arabic{bnomer}.\arabic{snomer}. }{#1}\par\medskip}
\newcommand{\theo}[2]{%
\refstepcounter{snomer}
\par\textbf{Theorem \arabic{bnomer}.\arabic{snomer}. }{#2} {\emph{#1}}\hspace{\fill}$\square$\par}
\newcommand{\mtheop}[2]{%
\refstepcounter{snomer}
\par\textbf{Theorem \arabic{bnomer}.\arabic{snomer}. }{\emph{#1}}
\par\textsc{Proof}. {#2}\hspace{\fill}$\square$\par}
\newcommand{\mcorop}[2]{%
\refstepcounter{snomer}
\par\textbf{Corollary \arabic{bnomer}.\arabic{snomer}. }{\emph{#1}}
\par\textsc{Proof}. {#2}\hspace{\fill}$\square$\par}
\newcommand{\mtheo}[1]{%
\refstepcounter{snomer}
\par\medskip\textbf{Theorem \arabic{bnomer}.\arabic{snomer}. }{\emph{#1}}\par\medskip}
\newcommand{\theobn}[1]{%
\par\medskip\textbf{Theorem. }{\emph{#1}}\par\medskip}
\newcommand{\theoc}[2]{%
\refstepcounter{snomer}
\par\medskip\textbf{Theorem \arabic{bnomer}.\arabic{snomer}. }{#1} {\emph{#2}}\par\medskip}
\newcommand{\mlemm}[1]{%
\refstepcounter{snomer}
\par\medskip\textbf{Lemma \arabic{bnomer}.\arabic{snomer}. }{\emph{#1}}\par\medskip}
\newcommand{\mprop}[1]{%
\refstepcounter{snomer}
\par\medskip\textbf{Proposition \arabic{bnomer}.\arabic{snomer}. }{\emph{#1}}\par\medskip}
\newcommand{\theobp}[2]{%
\refstepcounter{snomer}
\par\textbf{Theorem \arabic{bnomer}.\arabic{snomer}. }{#2} {\emph{#1}}\par}
\newcommand{\theop}[2]{%
\refstepcounter{snomer}
\par\textbf{Theorem \arabic{bnomer}.\arabic{snomer}. }{\emph{#1}}
\par\textsc{Proof}. {#2}\hspace{\fill}$\square$\par}
\newcommand{\theosp}[2]{%
\refstepcounter{snomer}
\par\textbf{Theorem \arabic{bnomer}.\arabic{snomer}. }{\emph{#1}}
\par\textsc{Sketch of the proof}. {#2}\hspace{\fill}$\square$\par}
\newcommand{\exam}[1]{%
\refstepcounter{snomer}
\par\medskip\textbf{Example \arabic{bnomer}.\arabic{snomer}. }{#1}\par\medskip}
\newcommand{\deno}[1]{%
\refstepcounter{snomer}
\par\textbf{Notation \arabic{bnomer}.\arabic{snomer}. }{#1}\par}
\newcommand{\lemm}[1]{%
\refstepcounter{snomer}
\par\textbf{Lemma \arabic{bnomer}.\arabic{snomer}. }{\emph{#1}}\hspace{\fill}$\square$\par}
\newcommand{\lemmp}[2]{%
\refstepcounter{snomer}
\par\medskip\textbf{Lemma \arabic{bnomer}.\arabic{snomer}. }{\emph{#1}}
\par\textsc{Proof}. {#2}\hspace{\fill}$\square$\par\medskip}
\newcommand{\coro}[1]{%
\refstepcounter{snomer}
\par\textbf{Corollary \arabic{bnomer}.\arabic{snomer}. }{\emph{#1}}\hspace{\fill}$\square$\par}
\newcommand{\mcoro}[1]{%
\refstepcounter{snomer}
\par\textbf{Corollary \arabic{bnomer}.\arabic{snomer}. }{\emph{#1}}\par\medskip}
\newcommand{\corop}[2]{%
\refstepcounter{snomer}
\par\textbf{Corollary \arabic{bnomer}.\arabic{snomer}. }{\emph{#1}}
\par\textsc{Proof}. {#2}\hspace{\fill}$\square$\par}
\newcommand{\nota}[1]{%
\refstepcounter{snomer}
\par\medskip\textbf{Remark \arabic{bnomer}.\arabic{snomer}. }{#1}\par\medskip}
\newcommand{\propp}[2]{%
\refstepcounter{snomer}
\par\medskip\textbf{Proposition \arabic{bnomer}.\arabic{snomer}. }{\emph{#1}}
\par\textsc{Proof}. {#2}\hspace{\fill}$\square$\par\medskip}
\newcommand{\hypo}[1]{%
\refstepcounter{snomer}
\par\medskip\textbf{Conjecture \arabic{bnomer}.\arabic{snomer}. }{\emph{#1}}\par\medskip}
\newcommand{\prop}[1]{%
\refstepcounter{snomer}
\par\textbf{Proposition \arabic{bnomer}.\arabic{snomer}. }{\emph{#1}}\hspace{\fill}$\square$\par}

\newcommand{\proof}[2]{%
\par\medskip\textsc{Proof{#1}}. \hspace{-0.2cm}{#2}\hspace{\fill}$\square$\par\medskip}

\makeatletter
\def\iddots{\mathinner{\mkern1mu\raise\p@
\vbox{\kern7\p@\hbox{.}}\mkern2mu
\raise4\p@\hbox{.}\mkern2mu\raise7\p@\hbox{.}\mkern1mu}}
\makeatother

\newcommand{\okr}[2]{%
\refstepcounter{snomer}
\par\medskip\textbf{{#1} \arabic{bnomer}.\arabic{snomer}. }{\emph{#2}}\par\medskip}

\newcommand{\Ind}[3]{%
\mathrm{Ind}_{#1}^{#2}{#3}}
\newcommand{\Res}[3]{%
\mathrm{Res}_{#1}^{#2}{#3}}
\newcommand{\epsi}{\varepsilon}
\newcommand{\tri}{\triangleleft}
\newcommand{\Supp}[1]{%
\mathrm{Supp}(#1)}
\newcommand{\SSu}[1]{%
\mathrm{SingSupp}(#1)}

\newcommand{\gee}{\geqslant}
\newcommand{\reg}{\mathrm{reg}}
\newcommand{\Dyn}{\mathrm{Dyn}}
\newcommand{\Ann}{\mathrm{Ann}\,}
\newcommand{\Cent}[1]{\mathbin\mathrm{Cent}({#1})}
\newcommand{\PCent}[1]{\mathbin\mathrm{PCent}({#1})}
\newcommand{\Irr}[1]{\mathbin\mathrm{Irr}({#1})}
\newcommand{\Exp}[1]{\mathbin\mathrm{Exp}({#1})}
\newcommand{\empr}[2]{[-{#1},{#1}]\times[-{#2},{#2}]}
\newcommand{\sreg}{\mathrm{sreg}}
\newcommand{\ilm}{\varinjlim}
\newcommand{\wdth}{\mathrm{wd}}
\newcommand{\plm}{\varprojlim}
\newcommand{\codim}{\mathrm{codim}\,}
\newcommand{\GKdim}{\mathrm{GKdim}\,}
\newcommand{\chara}{\mathrm{char}\,}
\newcommand{\rk}{\mathrm{rk}\,}
\newcommand{\chr}{\mathrm{ch}\,}
\newcommand{\Ker}{\mathrm{Ker}\,}
\newcommand{\id}{\mathrm{id}}
\newcommand{\Ad}{\mathrm{Ad}}
\newcommand{\Gh}{\mathrm{Gh}}
\newcommand{\col}{\mathrm{col}}
\newcommand{\row}{\mathrm{row}}
\newcommand{\high}{\mathrm{high}}
\newcommand{\low}{\mathrm{low}}
\newcommand{\pho}{\hphantom{\quad}\vphantom{\mid}}
\newcommand{\fho}[1]{\vphantom{\mid}\setbox0\hbox{00}\hbox to \wd0{\hss\ensuremath{#1}\hss}}
\newcommand{\wt}{\widetilde}
\newcommand{\wh}{\widehat}
\newcommand{\ad}[1]{\mathrm{ad}_{#1}}
\newcommand{\tr}{\mathrm{tr}\,}
\newcommand{\GL}{\mathrm{GL}}
\newcommand{\SL}{\mathrm{SL}}
\newcommand{\SO}{\mathrm{SO}}
\newcommand{\Or}{\mathrm{O}}
\newcommand{\Sp}{\mathrm{Sp}}
\newcommand{\Sa}{\mathrm{S}}
\newcommand{\Ua}{\mathrm{U}}
\newcommand{\Andre}{\mathrm{Andre}}
\newcommand{\Aord}{\mathrm{Aord}}
\newcommand{\Mat}{\mathrm{Mat}}
\newcommand{\Stab}{\mathrm{Stab}}
\newcommand{\htt}{\mathfrak{h}}
\newcommand{\spt}{\mathfrak{sp}}
\newcommand{\slt}{\mathfrak{sl}}
\newcommand{\sot}{\mathfrak{so}}

\newcommand{\vfi}{\varphi}
\newcommand{\aad}{\mathrm{ad}}
\newcommand{\vpi}{\varpi}
\newcommand{\teta}{\vartheta}
\newcommand{\Bfi}{\Phi}
\newcommand{\Fp}{\mathbb{F}}
\newcommand{\Rp}{\mathbb{R}}
\newcommand{\Zp}{\mathbb{Z}}
\newcommand{\Cp}{\mathbb{C}}
\newcommand{\Ap}{\mathbb{A}}
\newcommand{\Pp}{\mathbb{P}}
\newcommand{\Kp}{\mathbb{K}}
\newcommand{\Np}{\mathbb{N}}
\newcommand{\ut}{\mathfrak{u}}
\newcommand{\at}{\mathfrak{a}}
\newcommand{\glt}{\mathfrak{gl}}
\newcommand{\hei}{\mathfrak{hei}}
\newcommand{\nt}{\mathfrak{n}}
\newcommand{\kt}{\mathfrak{k}}
\newcommand{\mt}{\mathfrak{m}}
\newcommand{\rt}{\mathfrak{r}}
\newcommand{\rad}{\mathfrak{rad}}
\newcommand{\bt}{\mathfrak{b}}
\newcommand{\unt}{\underline{\mathfrak{n}}}
\newcommand{\gt}{\mathfrak{g}}
\newcommand{\vt}{\mathfrak{v}}
\newcommand{\pt}{\mathfrak{p}}
\newcommand{\Xt}{\mathfrak{X}}
\newcommand{\Po}{\mathcal{P}}
\newcommand{\PV}{\mathcal{PV}}
\newcommand{\Uo}{\EuScript{U}}
\newcommand{\Fo}{\EuScript{F}}
\newcommand{\Do}{\EuScript{D}}
\newcommand{\Eo}{\EuScript{E}}
\newcommand{\Jo}{\EuScript{J}}
\newcommand{\Iu}{\mathcal{I}}
\newcommand{\Mo}{\mathcal{M}}
\newcommand{\Nu}{\mathcal{N}}
\newcommand{\Ro}{\mathcal{R}}
\newcommand{\Co}{\mathcal{C}}
\newcommand{\Ko}{\mathcal{K}}
\newcommand{\So}{\mathcal{S}}
\newcommand{\Lo}{\mathcal{L}}
\newcommand{\Ou}{\mathcal{O}}
\newcommand{\Uu}{\mathcal{U}}
\newcommand{\Tu}{\mathcal{T}}
\newcommand{\Au}{\mathcal{A}}
\newcommand{\Vu}{\mathcal{V}}
\newcommand{\Du}{\mathcal{D}}
\newcommand{\Bu}{\mathcal{B}}
\newcommand{\Sy}{\mathcal{Z}}
\newcommand{\Sb}{\mathcal{F}}
\newcommand{\Gr}{\mathcal{G}}
\newcommand{\Xu}{\mathcal{X}}
\newcommand{\Et}{\mathcal{E}}
\newcommand{\Op}{\mathbb{O}}
\newcommand{\chv}{\mathrm{chv}}
\newcommand{\rtc}[1]{C_{#1}^{\mathrm{red}}}

\author{Mikhail Venchakov}
\date{}
\title{Rook placements and coadjoint orbits for maximal unipotent subgroups of finite symplectic groups}\maketitle
\begin{abstract} Let $U$ be a maximal unipotent subgroup in a symplectic group over a finite field of sufficiently large characteristic $p$. According to the Kirillov's orbit method, the coadjoint orbits of the group $U$ play the key role in the description of irreducible complex characters of $U$. Almost all important classes of orbits and characters studied to the moment can be uniformly described as the orbits and characters associated with so-called orthogonal rook placements. In the paper, we construct a semi-direct decomposition for the corresponding irreducible characters in the spirit of the Mackey little group method. As a corollary, we present an explicit formula for the character corresponding to an orbit of maximal possible dimension.

\end{abstract}

\sect{Introduction}

\let\thefootnote\relax\footnote{This work is an output of a research project implemented as part of the Basic Research Program at HSE University.}

Let $U$ be a unipotent algebraic group over a finite field $\Fp_q$ of sufficiently large characteristic $p$. The main tool in representation theory of $U$ is the orbit method created in 1962 by A.A. Kirillov, see~\cite{Kirillov62}, \cite{Kirillov04}, \cite{Kazhdan77}. The key idea of the orbit method says that the irreducible representations of $U$ are in one-to-one correspondence with the coadjoint orbits of this group. Namely, the group $U$ acts on its Lie algebra $\ut$ via the adjoint action; the dual action of $U$ on the dual space $\ut^*$ is called coadjoint. It turns out that there is a natural bijection between the set $\Irr{U}$ of all irreducible complex characters of $U$ and the set $\ut^*/U$ of coadjoint orbits, see Theorem~\ref{theo:orbit_method} below for the precise statement.

Let $U$ be a maximal unipotent subgroup (or, equivalently, a Sylow $p$-subgroup) in a simple classical group over $\Fp_q$. A complete description of all coadjoint orbits for the group $U$ is a wild problem, so a natural question is how to describe certain important classes of orbits and the corresponding irreducible characters. For the case $A_{n-1}$, orbits of maximal possible dimension were classified in the first Kirillov's work on the orbit method \cite{Kirillov62}. All of them are associated with the so-called Kostant cascade, which is in fact an orthogonal rook placement. For other classical root systems, orbits associated with the Kostant cascades also have maximal possible dimension, see \cite{Kostant12} and \cite{Kostant13}. For $A_{n-1}$, orbits of submaximal dimension were classified by A.N. Panov \cite{IgnatevPanov09}. Such orbits also correspond to orthogonal rook placements (modulo adding simple root covectors to the canonical forms on them). A classification of orbits of maximal dimension in type $C_n$ follows from C. Andre and A. Neto's paper \cite{AndreNeto06}. A description of orbits of maximal possible dimension for $B_n$ and $D_n$, as well as of orbits of submaximal dimension for all types $B_n$, $C_n$ and $D_n$, will be presented in \cite{IgnatevPetukhovVenchakov24}. All these orbits are also associated with orthogonal rook placements (modulo adding simple root covectors to the canonical forms on them).

The irreducible characters corresponding to the orbits of maximal dimension $M$ in type $A_{n-1}$ were computed by C. Andre in 2011 \cite{Andre01}. The characters of submaximal dimension $M-2$ for $A_{n-1}$ were calculated by M. Ignatev in \cite{Ignatev09}. A formula for the characters of maximal dimension in type $C_n$ can be obtained from the results of Andre and Neto \cite{AndreNeto06}. Characters of maximal dimension in the orthogonal case, as well as the characters of submaximal dimension for classical cases, will be computed in \cite{IgnatevPetukhovVenchakov24}. Andre's method is based on very nice stratification of $\ut^*$ established by him, whose strata are numerated by arbitrary (possibly, non-orthogonal) rook placements and maps from them to the set $\Fp_q^{\times}$ of nonzero numbers.

In this paper, we present another approach to calculating characters based on the Mackey method of little groups. This method reduces the calculation of irreducible characters of a semi-direct product group with abelian normal subgroup to the calculation of irreducible characters of certain subgroups called little groups, see Section~\ref{sect:sea_battle_pol} for the detail. It turned out that this method allows to compute characters of maximal and submaximal dimension for $A_{n-1}$, see \cite{Ignatev09}. Furthermore, in \cite{IgnatevVenchakov24} we applied this method to the calculation of characters of the next possible dimension $M-4$. Here we apply a modification of the Mackey method to the characters associated with arbitrary orthogonal rook placements in the symplectic case. As a corollary, we give an alternative proof of the formula for the characters of maximal dimension given in \cite{AndreNeto09}. We also give an independent classification of orbits of maximal dimension previously obtained in \cite{AndreNeto06}. 

The structure of the paper is as follows. In Section~\ref{main_definitions}, we briefly recall main definitions and facts about finite symplectic groups and its maximal unipotent subgroups. We also define rook placements in root systems and associated coadjoint orbits. In Section~\ref{sect:sea_battle_pol}, we recall the orbit method. As the main result, we present a decomposition for an arbitrary character corresponding to an orbit associated with orthogonal rook placement in the spirit of the Mackey method, see Theorem \ref{theo_char_decompose}. In particular, to construct the representation corresponding to the orbit of a linear form $f$, and to obtain a decomposition mentioned above, it is needed to present a polarization for $f$. In that section we also established such a polarization for the canonical form on an orbit associated with an orthogonal rook placement.

In the Section \ref{sect:char_max_dim}, we explore the characters of maximal dimension. At the first part, we provide a classification of orbits of maximal dimension in terms of so-called canonical forms, see Theorem~\ref{theo_orb_max_dim}. Next, we describe the support (i.e., the set of elements of group $U$ on which the character is not zero) of an arbitrary character of maximal dimension and calculate it's value. More precisely, we split the support into the union of conjugacy classes and calculate the value of our character on each class, see Theorem \ref{theo_char_value}. 

I thank my scientific advisor Mikhail Ignatev for constant attention to my work, and Alexey Petukhov for fruitful discussions.

\sect{Definitions and notations} \label{main_definitions}
The symplectic group has $C_n$ as it's corresponding root system. We will denote the root system of type $C_n$ by $\Phi$. As usual, we will identify it with the following subsets of $\Rp^n$:
\begin{center}
$C_n = \{\pm\epsi_i\pm\epsi_j, 1\leq i < j\leq
n\}\cup\{\pm2\epsi_i, 1\leq i\leq n\},$
\end{center}
where $\{\epsi_i\}_{i=1}^n$ is the standard basis in $\Rp^n$. Pick the set of simple roots $\Delta=\Delta(\Phi)$ as in
\cite{Bourbaki03}:

\begin{equation}
\Delta(C_n)=\{\epsi_i-\epsi_{i+1},1\leq i\leq n-1\}\cup\{2\epsi_n\}.
\end{equation}

The set of positive roots $\Phi^+\supset\Delta$ looks as follows:
\begin{equation}
C_n^+ = \{\epsi_i\pm\epsi_j, 1\leq i < j\leq n\}\cup\{2\epsi_i,
1\leq i\leq n\}.
\end{equation}

Denote by $\ut=\ut(\Phi)$ the subalgebra of $\mathfrak{gl}_m(\Fp_q)$
spanned by the vectors $e_{\alpha}$, $\alpha\in\Phi$, where
$$
\begin{array}{ll}
&e_{2\epsi_i}=e_{-i, i},\quad 1\leq i\leq n,\\
&e_{\epsi_i-\epsi_j}=e_{j,i}-e_{-i,-j},\quad 1\leq i<j\leq n,\\
&e_{\epsi_i+\epsi_j}=e_{-j,i}+e_{-i,j},\quad 1\leq i<j\leq n.\\
\end{array}
$$
Here, $m=2n$, and we numerate the rows and the columns of $m\times m$ matrix by the indices $$1,2,\ldots,n,-n,\ldots,-2,-1,$$ and we denote by $e_{a, b}$ the usual elementary matrix. Of course, it is a maximal nilpotent subalgebra in the corresponding symplectic algebra $\gt= \gt(\Phi)$. In particular, $\dim\ut = |\Phi^+|$. 

\exam{Here\label{ex:squares} we schematically drew the algebra $C_4$. The convention is as follows: given $1\leq j<i$, the square $(i,j)$ (respectively, $(-i,j)$ and $(-j,j)$) corresponds to the root $\epsi_j-\epsi_i$ (respectively, $\epsi_j+\epsi_i$ and $2\epsi_j$).
\begin{center}\small
$\mymatrix{
\lNote{1}\Note{1}\pho& \Note{2}\pho& \Note{3}\pho& \Note{4}\pho& \Note{-4}\pho& \Note{-3}\pho& \Note{-2}\pho& \Note{-1}\pho\\
\lNote{2}\Top{2pt}\Rt{2pt}& \pho& \pho& \pho& \pho& \pho& \pho& \pho\\
\lNote{3}\pho& \Top{2pt}\Rt{2pt}& \pho& \pho& \pho& \pho& \pho& \pho\\
\lNote{4}\pho& \pho& \Top{2pt}\Rt{2pt}& \pho& \pho& \pho& \pho& \pho\\
\lNote{-4}\pho& \pho& \pho& \Top{2pt}\Rt{2pt}\Bot{2pt}& \pho& \pho& \pho& \pho\\
\lNote{-3}\pho& \pho& \Bot{2pt}\Rt{2pt}& \pho& \Top{2pt}\Rt{2pt}& \pho& \pho& \pho\\
\lNote{-2}\pho& \Bot{2pt}\Rt{2pt}& \pho& \pho& \pho& \Top{2pt}\Rt{2pt}& \pho& \pho\\
\lNote{-1}\Bot{2pt}\Rt{2pt}& \pho& \pho& \pho& \pho& \pho& \Top{2pt}\Rt{2pt}& \pho\\
}$
\end{center}
}

Furthermore, we define the functions 
$$
\begin{array}{ll}
&\col\colon\Phi^+\to\{1,\ldots,n\}\colon\col(\epsi_i\pm\epsi_j)=\col(2\epsi_i)=i,\\
&\row\colon\Phi^+\to\{-n,\ldots,n\}\colon\row(\epsi_i\pm\epsi_j)=\mp j,\row(2\epsi_i)=-i.\\
\end{array}
$$
For arbitrary $-n+1\leq i\leq n-1$ and $1\leq j\leq n$, the sets
$$
\begin{array}{ll}
&R_i = R_i(\Phi) = \{\alpha\in\Phi^+\mid \row(\alpha)=i\},\\
&C_j = C_j(\Phi) = \{\alpha\in\Phi^+\mid \col(\alpha)=j\}\\
\end{array}
$$
are called the $i$th\emph{ row} and the $j$th\emph{ column}
$\Phi^+$ respectively. We introduce the mirror order on the 
set of indices $$1\prec2\prec\ldots\prec
n\prec-n\prec\ldots\prec-2\prec-1,$$ and the following total orders on $\Phi^+$:
\begin{equation}
\begin{array}{ll}
&\alpha\prec\beta~\stackrel{\mathrm{def}}{\iff}~
\col(\beta)\prec\col(\alpha)~\mbox{or}~
\col(\beta)=\col(\alpha),~\row(\beta)\prec\row(\alpha),\\
&\alpha\prec'\beta~\stackrel{\mathrm{def}}{\iff}~
\col(\beta)\prec\col(\alpha)~\mbox{or}~
\col(\beta)=\col(\alpha),~\row(\beta)\succ\row(\alpha).\label{formula_complete_orders}
\end{array}
\end{equation}
For example, for $\Phi=C_6$ we have
$\epsi_2-\epsi_4\succ\epsi_2\succ\epsi_2+\epsi_5\succ\epsi_3-\epsi_6$,
$2\epsi_2\succ'\epsi_2+\epsi_5\succ'\epsi_2-\epsi_4\succ'\epsi_3-\epsi_6$.

Now we are ready to define the main object of our interest. 

\defi{Let a subset $D=\{\beta_1,\ldots,\beta_t\}\subset\Phi^+$ consisting of pairwise orthogonal roots satisfy the following condition:
\begin{equation}
|D\cap R_i|\leq 1\mbox{ and }|D\cap C_j|\leq 1\mbox{ for all }i,
j.\label{formula_basic_supp}
\end{equation}
We call such a subset $D$ an \emph{orthogonal rook placement}. (Note than in \cite{AndreNeto06}, such a subset $D$ is called a basic subset of $\Phi$.)}

%\nota{Let $W=W(\Phi)$ be the Weyl group of the root system $\Phi$. To each orthogonal rook placement $D=\{\beta_1,\ldots,\beta_t\}$ one can attach the involution (i.e., the element of order 2)
%$\sigma=\sigma_D\in W$:
%\begin{equation}
%\sigma=r_{\beta_1}\ldots r_{\beta_t}.
%\label{formula_ortog_decompose}
%\end{equation}
%Here, for an arbitrary $\beta\in\Phi$, we denote by $r_{
%\beta}$ the reflection of $\Rp^n$ with respect to the hyperplane orthogonal to  $\beta$. We will call the set $D$ the \emph{support} of the  involution $\sigma$ and write $D=\Supp{\sigma}$.

%fix
%Note that, in the case of $A_n$ or $C_n$, each involution can be obtained in such a way. On the other hand, for $B_n$ or $D_n$, this is not true. For example, if $\Phi=B_2$ and $\sigma=r_{\epsi_1}\cdot r_{\epsi_2}$, then there is no orthogonal rook placement $D$ such that $D=\Supp{\sigma}$.

%Note that, in the case of $B_n$ or $D_n$, not each involution can be obtained in such a way. For example, if $\Phi=B_2$ and $\sigma=r_{\epsi_1}\cdot r_{\epsi_2}$, then there is no orthogonal rook placement $D$ such that $D=\Supp{\sigma}$.
%}

Denote maximal unipotent subgroup $\exp(\ut)$ in the corresponding symplectic finite group $G$ by~$U$. In the sequel, we will assume everywhere that $p=\chara{\Fp_q}\gee n$. Under this assumption, the map
$$
\exp\colon\ut\to U\colon x\mapsto \sum_{i=0}^{m-1}\frac{x^i}{i!}
$$
is well-defined and is in fact a bijection (and also an isomorphism of algebraic varieties over $\overline{\Fp_q}$);
we denote the inverse map by $\ln$. Furthermore, the 
Backer--Campbell--Hausdorff formula claims that, for a Lie subalgebra $\at\subset\ut$ and arbitrary $u, v\in\at$, one has
\begin{equation}
\exp(u)\exp(v)=\exp(u + v + \tau(u,
v)),\label{formula_Campbell_Hausdorf}
\end{equation} where
$\tau(u, v)\in[\at, \at]$ (here $[\at, \at]=\langle[x, y], x,
y\in\at\rangle_{\Fp_q}$).

The group $U$ acts on its Lie algebra $\ut$ by the adjoint action; the dual action of $U$ on the $\Fp_q$-dual space $\ut^*$ is called \emph{coadjoint}. Using the non-degenerate form $$\langle A, B\rangle=\mathrm{tr}(AB)$$ on $\mathfrak{gl}_m(\Fp_q)$, one can identify the dual space $\ut^*$ with the space $\ut^t$ (in this case, $e_{\alpha}^*=e_{\alpha}^t$ for any
$\alpha=2\epsi_i\in\Phi^+$ and $e^*_{\alpha}=e_{\alpha}^T/2$ otherwise). Under this identification, the coadjoint action has the following form:
$$
g.x = \mathrm{pr}(gxg^{-1}),~g\in U,~x\in\ut^*
$$
(here we denote the projection
$\mathfrak{gl}_m(\Fp_q)\to\ut^*$ along $\ut$ by $\mathrm{pr}$). The orbits of the coadjoint action play the key role in the classification of irreducible representations of group $U$.
\defi{Let $D=\{\beta_1,\ldots,\beta_t\}\subset\Phi^+$ be an orthogonal rook placement, and
$\xi=(\xi_{\beta})_{\beta\in D}$ be a set of nonzero constants from $\Fp_q$. Put
$$
f=f_{D, \xi} = \sum_{\beta\in D}\xi_{\beta}e_{\beta}^*\in\ut^*. $$ We say that the coadjoint orbit $\Omega=\Omega_{D,
\xi}\subset\ut^*$ of the linear form $f$ is \emph{associated} with $D$, and $f$ is called the
\emph{canonical form} on $\Omega$.}

\defi{For an arbitrary Lie algebra $\gt$ and a linear form $f\in\gt^*$ let us consider a set $\pt$ which satisfies the following conditions.
\begin{enumerate}
    \item For an arbitrary $x,y\in\pt$, $f([x,y])=0$.
    \item The set $\pt$ is maximal subspace in $\ut$ satisfying property 1.
    \item The set $\pt$ is a subalgebra.
\end{enumerate}
Such a subspace is called a \emph{polarization} for the form $f$. 

Note that for any linear form a polarization always exists \cite{Vergne}. 
Also, it is easy to see that the dimension of the coadjoint orbit of $f$ is connected with the dimension of a polarization of $f$ by the simple formula $$\dim\Omega_f=2\cdot\codim\pt$$ (see, e.g., \cite[page
117]{Srinivasan}).}

\newpage\sect{Mackey type decomposition for characters\break associated with orthogonal rook placements} \label{sect:sea_battle_pol}
First, I want to recall the main idea of the orbit method created by A.A. Kirillov \cite{Kirillov62}, \cite{Kirillov04} and adapted by D. Kazhdan \cite{Kazhdan77} to finite groups (see also \cite{Kirillov95}), and provide a construction for orbits associated with orthogonal rook placements in symplectic groups. According to the orbit method, the construction of a polarization is a necessary step to obtain the irreducible representation corresponding to the orbit of a given linear form.

Specifically, pick and fix $\theta\colon\Fp_q\to\Cp^{\times}=\Cp\setminus\{0\}$, an arbitrary nontrivial character of the field $\Fp_q$ (i.e., a nontrivial homomorphism from the additive group of this field to the multiplicative group of complex numbers). Recall that we called the map inverse to $\exp\colon\ut\to U$ by $\ln\colon
U\to\ut$.
Let $\Omega\subset\ut^*$ be an arbitrary coadjoint orbit.
Consider the following function $\chi=\chi_{\Omega}\colon U\to\Cp$:
\begin{equation}
\label{formula_char_orbit}
\chi(g)=q^{-\frac{1}{2}\dim\Omega}\cdot\sum_{f\in\Omega}\theta(f(\ln(g))),\quad
g\in U
\end{equation}
(over $\overline{\Fp}_q$, coadjoint orbits are affine varieties of even dimension in $\ut^*$; furthermore,
$q^{\dim\Omega}=|\Omega|$).
The main idea of the orbit method can be expressed by the following theorem. \theo{The map\label{theo:orbit_method}
$\Omega\mapsto\chi_{\Omega}$ establishes a one-to-one correspondence between the set of coadjoint orbits of the group $U$ and the set $\Irr{U}$ of finite-dimensional irreducible complex characters of the group $U$. Furthermore\textup, the complex dimension of the representation corresponding to $\Omega$ equals $$q^{\frac{1}{2}\dim\Omega}=\sqrt{|\Omega|}.$$}{\cite[Proposition
2]{Kazhdan77}} \label{theo_char_orbit} Moreover, the construction of the representation corresponding to a given orbit $\Omega$ is as follows. One should pick a linear form $f\in\Omega$ and a polarization $\pt$ for $f$. Next, one should consider the following one-dimensional representation of $\exp(\pt)$:
$$\vfi_f(h)=\theta(f(\ln(h))),~h\in\exp(\pt).
\label{formula_irr_char}$$
The required representation now has the form $\Ind{\exp(\pt)}{U}(\vfi_f)$.

Now, we aim to provide a construction of a polarization for an arbitrary linear form $f=f_{D,\xi}$ associated with an orthogonal rook placement. Let us introduce some more notations.

\defi{Let $\beta\in\Phi^+$. Roots $\alpha$, $\gamma\in\Phi^+$ are called
$\beta$-\emph{singular} if $\alpha+\gamma=\beta$. The set of all $\beta$-singular roots is denoted by $S(\beta).$}

It is easy to see that the singular roots look as follows:
\begin{equation*}
\begin{split}
S(\epsi_i-\epsi_j)=&\bigcup_{l=i+1}^{j-1}\{\epsi_i-\epsi_l,\epsi_l-\epsi_j\},~1\leq i < j\leq n,\\
S(2\epsi_i) = &\bigcup_{l=i+1}^{n}\{\epsi_i-\epsi_l,\epsi_i+\epsi_l\},~1\leq i\leq n,\\
%&S(2\epsi_i) = \bigcup_{j=i+1}^{n}\{\epsi_i-\epsi_j,\epsi_i+\epsi_j\},~1\leq i\leq n,\\
S(\epsi_i+\epsi_j)=&\bigcup_{l=i+1}^{j-1}\{\epsi_i-\epsi_l,\epsi_l+\epsi_j\}
\cup\bigcup_{l=j+1}^n\{\epsi_i-\epsi_l,\epsi_j+\epsi_l\}\\
\cup&\bigcup_{l=j+1}^n\{\epsi_j+\epsi_l,\epsi_j-\epsi_l\}\cup S_{ij}
,~1\leq i<j\leq
n,\\
\end{split}\label{formula_sing_roots}
\end{equation*}
where $S_{ij}=\{\epsi_i-\epsi_j,2\epsi_j\}.$

\nota{Note that, given $\beta\in\Phi^+$, %$\beta\neq2\epsi_i$,
the column of $\beta$ contains exactly one of the root from each pair of $\beta$-singular roots whose sum is $\beta$, except for the case $\beta=2\epsi_i$. Precisely, put
\begin{equation*}
S^+(\beta)=
\begin{cases}
\{\epsi_i+\epsi_j,~i<k\leqslant n\},&\text{ if }\beta=2\epsi_i,\\
S(\beta)\cap C_{\text{col}(\beta)}&\text{ otherwise},
\end{cases}
\end{equation*}
so that $S(\beta)=S^+(\beta)\cup S^-(\beta)$\label{nota_one_row_col}, where $S^-(\beta)=S(\beta)\backslash S^+(\beta)$.}

\deno{Let $j_1<\cdots<j_t$ be the numbers of columns of roots from an orthogonal rook placement $D$. Put $j_0=0$ and define $M=M_D=\bigcup_{i=0}^tM_i$, where $M_0=\varnothing$ and, for any $i=1,\ldots,t$, inductively,
\begin{equation*}
    M_{j_i}=\{\gamma\in S^-(\beta)\mid\beta\in D\cap C_{j_i}\text{ and }\gamma,\beta-\gamma\notin\bigcup_{l=0}^{i-1}M_{j_l}\}.
\end{equation*}}
The next proposition provides a construction of a polarization for $f_{D,\xi}$ \cite[Theorem 1.1]{Ignatev2009}.
\prop{Let us consider $f=f_{D,\xi}$ and put $P=\Phi^+\backslash M$\textup{,} then 
\begin{equation*}
    \pt=\pt_{D,\xi}=\sum_{\alpha\in P}ke_{\alpha}
\end{equation*}
is a polarization for the form $f$.}

Above we observed a connection between coadjoint orbits and irreducible representations. The orbit method also provides us a formula \eqref{formula_char_orbit} for the irreducible character corresponding to a coadjoint orbit, but this formula is not convinient for calculations in case we want to obtain an explicit formula for character, so we have to use some special methods. More precisely, in the papers \cite{Ignatev09} and \cite{IgnatevVenchakov24} we used the Mackey method of little groups to calculate the value of subregular characters and characters of depth 2 for $U$ of type $A_n$. Furthermore, in the paper \cite{IgnatevVenchakov25} for cases when group $U$ has type $B_n$ and $D_n$, it was used to obtain so-called Mackey decomposition for characters corresponding to the orbits associated with orthogonal rook placements. Here we modify the Mackey method for type $C_n$. The Mackey method states the following. 

Let $\Gr$ be an arbitrary finite group, $A$, $B$ be its subgroups such that $\Gr$ is the semi-direct product of $A$ and $B$, i.e., $\Gr=AB$, $A\triangleleft\Gr$ and
$A\cap B=\{1\}$ (we denote $\Gr=A\rtimes B$). Assume additionally that $A$ is an abelian group and that $\psi$ is its irreducible character. The \emph{centralizer}
of $\psi$ in the group $B$ (or the \emph{little group}) is the subgroup of the form $$B^{\psi}=\{b\in
B\mid\psi\circ\tau_b=\psi\},$$ where $\tau_b\colon A\to A.~a\mapsto bab^{-1}.$

An arbitrary element $g\in\Gr$ can be uniquely represented as
$g=ab$, where $a\in A, b\in B$; this defines the maps
$\pi^{\Gr}_{A}\colon\Gr\to A\colon g\mapsto a,~\pi^{\Gr}_{B}\colon\Gr\to B\colon g\mapsto b,~g\in\Gr.$ Note that $\pi^{\Gr}_B$ is a group homomorphism, while $\pi^{\Gr}_A$, in general, is not. For an arbitrary subgroup $\wt B$ in~$B$ and arbitrary characters
$\psi\colon A\to\Cp$ and $\eta\colon\wt B\to\Cp$, one can define the character $\psi_0\eta_0$ of the group
$A\rtimes\wt B=A\wt B$, where
$$
\psi_0=\psi\circ\pi^{A\rtimes\wt
B}_{A},~\eta_0=\eta\circ\pi^{A\rtimes\wt B}_{\wt B}.
$$
Denote by $\Irr{\Gr}$ the set of irreducible characters of $\Gr$. The basic idea of the Mackey method, which allows us
to reduce the study of representations of the group $\Gr$ to the study of representations of
its subgroups $A$ and $B$, can be formulated as follows. \theo{Let $\Gr=A\rtimes B$
be a finite group with $A$ being an abelian group. Then each irreducible character $\chi$ of the group $\Gr$ has the form
\begin{equation}
\chi=\Ind{A\rtimes B^{\psi}}{G}{\psi_0\eta_0}
\label{formula_semi_direct}
\end{equation}
for certain $\psi\in\Irr{A}$\textup, $\eta\in\Irr{B}$. Conversely\textup, every function of the form \textup{(\ref{formula_semi_direct})} is an irreducible character of $\Gr$.}{\cite[Proposition
1.2]{Lehrer}}\label{theo_semi_direct}\newpage
%\postdiplaypenalty{10000}

However, in the case $C_n$ we have a problem with finding appropriate groups $A$ and $B$ such that \break $U=A\rtimes B$ and $A$ is abelian, which will be described later. So we can not use the Mackey method directly, but we are going to provide a decomposition of an arbitrary character associated with orthogonal rook placement in the spirit of the Mackey method with minor modifications. 

We will consider two separated cases. For an orbit $\Omega_{D,\xi}$ associated with an orthogonal rook placement $D=\{\beta_1\succ\ldots\succ\beta_t\}$ in $\Phi^+$, let us denote by $\ut_1$ and $\vt$ the subalgebras of $\ut$ of the following form:
\begin{equation*}
    \ut_1=\bigoplus_{\alpha\in
C'}\Fp_qe_{\alpha}
,~\vt=\bigoplus_{\alpha\in\Phi^+\setminus
C_1}\Fp_qe_{\alpha},
\end{equation*}
where $C'=C_1\backslash\{\epsi_1-\epsi_2,\epsi_1-\epsi_3,\ldots,\epsi_1-\epsi_n\}$ if $D\cap C_1=2\epsi_1$, and $C'=C_1$ otherwise. Note that in the case when $D\cap C_1=2\epsi_1$, the algebra $\ut_1$ is obviously not abelian, and $\ut\not=\ut_1\oplus\vt$ otherwise. Denote $U_1=\exp(\ut_1)$,
$V=\exp(\vt)$. Also define a function $\psi=\psi_{D,\xi}$ on the group $U_1$ by
\begin{equation}
\psi(x)=\theta(f(\ln(x)))=\begin{cases}
\theta(\xi_{\beta_1}
x_{i,1}),&\mbox{if }\col(\beta_1)=1,\row(\beta_1)=i,\\
1,&\mbox{if }\col(\beta_1)>1.
\end{cases}\label{formula_char_psi}
\end{equation}
As usual, we denote the $(i,j)$th element of a matrix $x$ by $x_{ij}$.

On the other hand, consider the following subgroup of $V$:
$$
\begin{array}{ll}
&V'=\prod_{\alpha\in \Phi_1^+}X_{\alpha},\mbox{ where}\\
&\Phi_1^+=\left\{\begin{array}{ll} \Phi^+\setminus C_1,&\mbox{if
}\col(\beta_1)>1,\\
\Phi^+\setminus(C_1\cup S^-(\beta_1)),&\mbox{if }\col(\beta_1)=1,
\end{array}\right.\\
&X_{\alpha}=\{x_{\alpha}(t),t\in\Fp_q\},~x_{\alpha}(t)=\exp(te_{\alpha}).
\end{array}
$$
In particular, if $\col(\beta_1)>1$ or $D\cap C_1=2\epsi_1$, then $V'=V$. It is clear that $V'=\exp(\vt')$, where $\vt'=\bigoplus_{\alpha\in \Phi_1^+}\Fp_qe_{\alpha}$.

\lemmp{For an arbitrary orthogonal rook placement
$D\subset\Phi^+$ and an arbitrary $\xi$\textup, the subgroup $V'$
coincides with the centralizer of the character $\psi$ in the subgroup $V$\textup, i.e.\textup, $V'=V^{\psi}$.\label{lemm_stab_psi}}{
It is clear if $\col(\beta_1)>1$, because in this case $V'=V^{\psi}=V$. Let $$\col(\beta_1)=1, \row(\beta_1)=i.$$ Here and next, we will write $x=\exp(y)\in V$ and $h=\exp(z)\in U_1$, where
$$
y=\sum_{\alpha\in\Phi^+\setminus
C_1}y_{\alpha}e_{\alpha}\in\vt,~z=\sum_{\gamma\in
C'}z_{\gamma}e_{\gamma}\in\ut_1,~y_{\alpha},z_{\gamma}\in\Fp_q.
$$
Since $h=1_m+z$, one has
$$
\begin{array}{ll}
xhx^{-1}&=\exp(y)(1_m+z)(\exp(y))^{-1}\\
&=1_m+\exp(y)z(\exp(y))^{-1}.
\end{array}
$$
But $\exp(y)z(\exp(y))^{-1}=(\exp\ad{y})(z)$, so we obtain
$$
xhx^{-1}=1_m+(\exp\ad{y})(z).
$$
In particular, by formula (\ref{formula_char_psi}), $$\psi(xhx^{-1})=\theta(\xi_{\beta_1}((\exp\ad{y})(z))_{i,1}).$$
Note that $V'=V$ in case $\beta_1=2\epsi_1$, so $V^{\psi}\subset V'$. It means that we need only to prove that $V'\subset V^{\psi}$. It is clear whenever $\col(\beta_1)=1$. Indeed, by definition, $V'$ is the product of the subgroups $X_{\alpha}$ for $\alpha\in\Phi_1^+$, but this subgroups centralize the character $\psi$. Hence, it is left to prove that $V^{\psi}\subset V'$ in the case when $\beta_1\ne 2\epsi_1$.

Assume that $y\in\vt\backslash\vt'$, we claim that there exists $z\in\ut_1$ such that $z_{i,1}\ne((\exp\ad y)(z))_{i,1}$. More precisely, let us rewrite $y$ in the following form:
$$
y=\sum_{\alpha\in \Phi_1^+}y_{\alpha}e_{\alpha}+\sum_{\delta\in
S^-(\beta_1)}y_{\delta}e_{\delta}\in\vt\setminus\vt',~
y_{\alpha},~y_{\delta}\in \Fp_q.
$$
Let $\delta_0$ be the smallest root with respect to the order $\prec'$ (see
(\ref{formula_complete_orders})) from all $\delta\in S^-(\beta_1)$
for which $y_{\delta}\neq0$. There exists $\gamma_0\in\Phi^+$ such that $\gamma_0+\delta_0=\beta_1$ (furthermore, $\gamma_0\in C_1$). Consider the element $z=e_{\gamma_0}\in\ut_1$. By definition, $z_{i,1}=0$. At the same time, $$
(\exp\ad{y})(z)=z+\ad{y}z+\frac{1}{2}\ad{y}^2z+\ldots.
$$
At first, let us assume that $\delta_0\ne2\epsi_i$ (the case $\delta_0=2\epsi_i$ is little technically different, we will consider it later). It is easy to see that $\ad{y}z=c_0y_{\delta_0}e_{\beta_1}$ where $c_0\ne0$, so $\ad{y}^2z=0$ because $\beta_1+\delta\notin\Phi^+$ for all $\delta\in S^-(\beta_1)$. Hence, $$((\exp\ad{y})(z))_{i,1}=c_0y_{\delta_0}\ne 0.$$

Now, let us move to the case when $\delta_0=2\epsi_i$. Define the set $S^*_i=S^*(\beta_1,\delta_0)$ to be the set of roots from $S^-(\beta_1)$ which are bigger than $\delta_0$ in sense of the order $\prec$. The roots from $S^*_i$ are, in fact, roots from the column $C_i$ without the root $2\epsi_i$. Here we see that $\gamma_0=\epsi_1-\epsi_i$. The obvious calculations shows us that $$\ad{y}z=c_0y_{\delta_0}e_{\beta_1}+\sum\limits_{S^*_i}y_{\delta}e_{\delta+\gamma_0}.$$
Note that all the roots $\delta+\gamma_0$ where $\delta\in S^*_i$, belong to $\Phi^+$ and even to $C_1$, more precisely, these are roots from $S^*_i$, whose columns indexes $i$ are replaced by 1. Since $2\epsi_i-\delta+\gamma_0+\delta\in\Phi^+$ (moreover, this root is $\beta_1$) and the root $2\epsi_i-\delta+\gamma_0$ is the only positive root which gives a positive root being added to $\delta$, we obtain that $$\ad{y}^2z=\sum\limits_{S^*_i}y_{\delta}y_{2\epsi_i-\delta}[e_{\delta},e_{2\epsi_i-\delta+\gamma_0}]=\sum\limits_{S^*_i}y_{\delta}y_{2\epsi_i-\delta}c_\delta e_{\beta_1}.$$
Here $c_{\delta}=1$ if $\row(\delta)>0$ and $c_{\delta}=-1$ otherwise. Hence, this sum is zero and we obtain that $$((\exp\ad{y})(z))_{i,1}=c_0y_{\delta_0}\ne 0$$
for some nonzero constant $c_0$. That completes the proof.
}

So, we have calculated the centralizer of the character $\psi$ in the subgroup $V$. Let us study its structure in more detail. Denote
$$
\Phi_2^+=\left\{\begin{array}{ll} \Phi_1^+\setminus
S^-(\epsi_1\pm\epsi_i),&\mbox{if
}\beta_1=\epsi_1\mp\epsi_i,~2\leq i\leq n,\\
\Phi_1^+&\mbox{otherwise},
\end{array}\right.
$$
and put
\begin{equation}
\vt_1=\bigoplus_{\alpha\in\Phi_1^+\setminus\Phi_2^+}\Fp_qe_{\alpha},~
\wt\ut=\bigoplus_{\alpha\in\Phi_2^+}\Fp_qe_{\alpha}.\label{formula_v_1}
\end{equation} It is easy to check that $\vt'=\vt_1\oplus\wt\ut$, where $\vt_1$ is an ideal of $\vt'$. Hence, $V'=V_1\rtimes\wt U$, where
$$V_1=\exp(\vt_1)=\prod_{\alpha\in\Phi_1^+\setminus\Phi_2^+}X_{\alpha},~\wt U=\exp(\wt\ut)=\prod_{\alpha\in\Phi_2^+}X_{\alpha}$$ (if
$\Phi_1^+=\Phi_2^+$ then $V_1=\{1_m\}$ and $\wt U=V'$).

Moreover, we notice that $\wt U$ is isomorphic to maximal unipotent subgroup of the symplectic group of smaller rank. More precisely, set
$$
\wt\Phi=\begin{cases}
C_{n-1},&\mbox{if }\col(\beta_1)>1\text{ or }D\cap C_1=2\epsi_1,\\
C_{n-2}&\mbox{otherwise},\\
\end{cases}.
$$
then one can easily see that $\wt\ut\cong\ut(\wt\Phi)$. There also exists the obvious one-to-one correspondence between root systems $\Phi_2^+$ and $\wt\Phi$, we will denote it by $\pi$. Slightly abusing the notations we will denote the corresponding isomorphism between algebras $\wt\ut$ and $\ut(\wt\Phi)$ by the same letter $\pi$.

Now, finally, we can present a decomposition for the irreducible character associated with an orthogonal rook placement. Recall that we have defined the function $\psi$ on the group $U_1$ and found
its centralizer $V^{\psi}=V'=V_1\rtimes\wt U$. It can be directly checked that $\psi\circ\pi_{U_1}^{U_1\rtimes V'}$ is a homomorphism, so it is a 1-dimensional representation of the group $U_1\rtimes V'$ and its character. In this case, $\wt\ut=\mathrm{Lie}(\wt U)\cong\ut(\wt\Phi)$, so
the irreducible characters of the group $\wt U$ are in
one-to-one correspondence with the coadjoint orbits in $\wt\ut^*$ (cf. (\ref{formula_char_orbit})).

As above, let  $D=\{\beta_1\succ\ldots\succ\beta_t\}$~ be
an orthogonal rook placement in $\Phi^+$, $\xi=(\xi_{\beta})_{\beta\in D}$,
$\Omega=\Omega_{D,\xi}$ be the associated orbit, and $\chi$
be the corresponding irreducible character of the group $U$. We denote $\wt D
=\pi(D\setminus C_1)\subset\wt\Phi^+$,
$\wt\Omega=\Omega_{\wt D,\wt\xi}$, and let $\wt\chi$ be the corresponding irreducible character of the group $\wt U$. Here
$$
\wt\xi=\begin{cases}
\xi\in(\Fp_q^{\times})^t,&\mbox{if }\col(\beta_1)>1,\\
(\xi_{\beta_2},\ldots,\xi_{\beta_t})\in(\Fp_q^{\times})^{t-1},&\mbox{if
}\col(\beta_1)=1.
\end{cases}.
$$
Denote by $\eta$ the homomorphism defining the representation $V_{\wt\chi}$ corresponding to the character $\wt\chi$. Then $\eta$ is a homomorphism, so $\eta\circ \pi_{\wt U}^{V'}\circ\pi_{V'}^{U_1\rtimes V'}$ is a homomorphism defining a representation of $U_1\rtimes V'$ with the character $\wt\chi\circ \pi_{\wt U}^{V'}\circ\pi_{V'}^{U_1\rtimes V'}$. Thus, the multiplication of the characters $(\wt\chi\circ \pi_{\wt U}^{V'}\circ\pi_{V'}^{U_1\rtimes V'})\cdot(\psi\circ\pi_{U_1}^{U_1\rtimes V'})$ of the group $U_1\rtimes V'$ will be a character of this group as a character of the tensor product of the corresponding representations.

\theop{Let $\pi^{U_1\rtimes V'}_{\wt U}=\pi^{V'}_{\wt
U}\circ\pi^{U_1\rtimes V'}_{V'}$. Then \begin{equation}
\chi=\Ind{U_1\rtimes V'}{U}{((\psi\circ\pi^{U_1\rtimes
V'}_{U_1})\cdot(\wt\chi\circ\pi^{U_1\rtimes V'}_{\wt
U}))}.\label{formula_char_decompose}
\end{equation}\label{theo_char_decompose}}{Put $P=\exp(\pt)\subset U, \wt
P=\exp(\wt\pt)\subset\wt U$, then $P=U_1\rtimes(V_1\rtimes\wt
P)$. Let us denote the function on the right-hand side of
(\ref{formula_char_decompose}) by $\eta$. By Theorem
\ref{theo_char_orbit},
$$
\chi=\Ind{P}{U}{(\theta\circ f\circ\ln)},\quad \wt\chi=\Ind{\wt
P}{\wt U}{(\theta\circ f\circ\ln)}.
$$
We will break down further reasoning into a series of steps.

1. We claim that if $y\in P$ then $f(\ln\pi^{P}_{V_1\rtimes\wt
P}(y))=f(\ln\pi^P_{\wt P}(y))$. Indeed, it is evident that
$$
\pi^{P}_{V_1\rtimes\wt P}(y)=(\pi^{V_1\rtimes\wt
P}_{V_1}\circ\pi^{P}_{V_1\rtimes\wt P})(y)\cdot(\pi^{V_1\rtimes\wt
P}_{\wt P}\circ\pi^{P}_{V_1\rtimes\wt
P})(y)=\pi^P_{V_1}(y)\cdot\pi^P_{\wt P}(y).
$$
Let $u=\ln\pi^P_{V_1}(y)\in\pt$, $v=\ln\pi^P_{\wt P}(y)\in\pt$, then, from (\ref{formula_Campbell_Hausdorf}), one has
$$\pi^{P}_{V_1\rtimes\wt P}(y)=\exp(u)\exp(v)=\exp(u+v+\tau),
\quad\mbox{where }\tau\in[\pt,\pt],$$ so
$\ln\pi^{P}_{V_1\rtimes\wt P}(y)=\ln\pi^P_{V_1}(y)+\ln\pi^P_{\wt
P}(y)+\tau$. But $\pt$ is a polarization for $f$, so $f(\tau)=0$. At the same time, $\ln\pi^P_{V_1}(y)\in\vt_1$ and, clearly,
$f{\mid}_{\vt_1}\equiv0$ (see the definition of $\vt_1$ in
(\ref{formula_v_1})). It follows that
$$
f(\ln\pi^{P}_{V_1\rtimes\wt
P}(y))=f(\ln\pi^P_{V_1}(y))+f(\ln\pi^P_{\wt
P}(y))+f(\tau)=f(\ln\pi^P_{\wt P}(y)).
$$

2. We claim that if $y\in P$, then $f(\ln y-\ln\pi^P_{\wt
P}(y)-\ln\pi^P_{U_1}(y))=0$. The proof is similar to the previous step. Indeed, if $u=\ln\pi^P_{V_1\rtimes \wt P}(y)\in\pt$,
$v=\ln\pi^P_{U_1}(y)\in\pt$ then, according to
(\ref{formula_Campbell_Hausdorf}),
$$
y=\pi^P_{U_1}(y)\cdot\pi^P_{V_1\rtimes\wt
P}(y)=\exp(u)\exp(v)=\exp(u+v+\tau),\quad\mbox{where }\tau\in[\pt,\pt].
$$
Hence, $\ln y=\ln\pi^P_{U_1}(y)+\ln\pi^P_{V_1\rtimes\wt P}+\tau$.
But $\pt$ is a polarization for $f$, so $f(\tau)=0$ and, consequently,
$f(\ln\pi^{P}_{V_1\rtimes\wt P}(y))=f(\ln\pi^P_{\wt P}(y))$.

3. The next claim is as follows. Let $\Gr= A\rtimes B$ be a finite group, $C$ be a subgroup in $B$ and $\lambda$ be a complex representation of the group $C$. Then
$\Ind{A\rtimes C}{G}{(\lambda\circ\pi^{A\rtimes
C}_C)}=(\Ind{C}{B}{\lambda})\circ\pi^{\Gr}_B$. This is proved in \cite[Proposition
1.2]{Lehrer}.

4. Let us apply step 3 to the case $\Gr=V_1\rtimes\wt U$, $A=V_1$, $B=\wt
U$, $C=\wt P$, $\lambda=\theta\circ f\circ\ln$. We obtain that
$$
\Ind{\wt P}{\wt U}{(\theta\circ f\circ\ln)}\circ\pi^{V_1\rtimes\wt
U}_{\wt U}=\Ind{V_1\rtimes\wt P}{V_1\rtimes\wt U}{(\theta\circ
f\circ\ln\circ\pi^{V_1\rtimes\wt P}_{\wt P})}.
$$
Now we apply step 3 to the case of $\Gr=U_1\rtimes V'$,
$A=U_1$, $B=V'=V_1\rtimes\wt U$, $C=V_1\rtimes\wt P$,
$\lambda=\theta\circ f\circ\ln\circ\pi^{V_1\rtimes\wt P}_{\wt P}$.
We have
$$
\begin{array}{ll}
\Ind{V_1\rtimes\wt P}{V_1\rtimes\wt U}{(\theta\circ
f\circ\ln\circ\pi^{V_1\rtimes\wt P}_{\wt P})}\circ\pi^{U_1\rtimes
V'}_{V_1\rtimes\wt U}&=\Ind{U_1\rtimes(V_1\rtimes\wt P)}{U_1\rtimes
V'}{(\theta\circ f\circ\ln\circ\pi^{V_1\rtimes\wt P}_{\wt
P}\circ\pi^{U_1\rtimes (V_1\rtimes \wt P)}_{V_1\rtimes\wt P})}=\\
&=\Ind{P}{U_1\rtimes V'}{(\theta\circ f\circ\ln\circ\pi^P_{\wt P})}.
\end{array}
$$
Hence, keeping in mind that $V'=V_1\rtimes\wt U$, we obtain
$$
\begin{array}{ll}
\Ind{\wt P}{\wt U}{(\theta\circ f\circ\ln)}\circ\pi^{U_1\rtimes
V'}_{\wt U}&=(\Ind{\wt P}{\wt U}{(\theta\circ
f\circ\ln)})\circ\pi^{V_1\rtimes\wt U}_{\wt U}\circ\pi^{U_1\rtimes
V'}_{V_1\rtimes \wt U}=\\
&=\Ind{P}{U_1\rtimes V'}{(\theta\circ f\circ\ln\circ\pi^P_{\wt P})}.
\end{array}
$$

5. One can rewrite the right-hand side of (\ref{formula_char_decompose}) taking into account
step 4 in the following way:
$$
\begin{array}{ll}
\eta&=\Ind{U_1\rtimes V'}{U}{((\psi\circ\pi^{U_1\rtimes
V'}_{U_1})\cdot(\wt\chi\circ\pi^{U_1\rtimes V'}_{\wt
U}))}=\\
&=\Ind{U_1\rtimes V'}{U}{((\theta\circ f\circ\ln\circ\pi^{U_1\rtimes
V'}_{U_1})\cdot(\Ind{\wt P}{\wt U}{(\theta\circ
f\circ\ln)}\circ\pi^{U_1\rtimes V'}_{\wt U}))}=\\
&=\Ind{U_1\rtimes V'}{U}{((\theta\circ f\circ\ln\circ\pi^{U_1\rtimes
V'}_{U_1})\cdot\Ind{P}{U_1\rtimes V'}{(\theta\circ
f\circ\ln\circ\pi^P_{\wt P})})}.
\end{array}
$$

6. We claim that the following equality holds:
$$
\Ind{P}{U_1\rtimes V'}{(\theta\circ f\circ\ln)}=(\theta\circ
f\circ\ln\circ\pi^{U_1\rtimes V'}_{U_1})\cdot\Ind{P}{U_1\rtimes
V'}{(\theta\circ f\circ\ln\circ\pi^P_{\wt P})}.
$$
Indeed, let $H\subset U_1\rtimes V'$ be an arbitrary complete
representative system of $(U_1\rtimes V')/P$ (in fact, one can choose $H$ to be a subset of $V'$, because $U_1\subset P$). Then, for any $x\in U_1\rtimes V'$,
$$
\begin{array}{ll}
&((\theta\circ f\circ\ln\circ\pi^{U_1\rtimes
V'}_{U_1})\cdot\Ind{P}{U_1\rtimes V'}{(\theta\circ
f\circ\ln\circ\pi^P_{\wt P})})(x)=\\
&=\theta(f(\ln\pi^{U_1\rtimes V'}_{U_1}(x)))\cdot\sum_{h\in H\mid
h^{-1}xh\in
P}\theta(f(\ln\pi^{P}_{\wt P}(h^{-1}xh)))=\\
&=\sum_{h\in H\mid h^{-1}xh\in P}\theta(f(\ln\pi^{U_1\rtimes
V'}_{U_1}(x)+\ln\pi^P_{\wt P}(h^{-1}xh))).
\end{array}
$$
But $h\in B$ implies that $h^{-1}\pi^{\Gr}_{A}(x)h=\pi^{\Gr}_A(h^{-1}xh)$
for an arbitrary finite group $\Gr=A\rtimes B$ and arbitrary
$x\in G$, $h\in B$. In our situation, this means that
$$\psi(X)=f(\ln\pi^{U_1\rtimes V'}_{U_1}(x))=f(\ln(h^{-1}\pi^{U_1\rtimes
V'}_{U_1}(x)h))=f(\ln\pi^{U_1\rtimes V'}_{U_1}(h^{-1}xh)),$$
because $h\in H\subset V' = V^{\psi}$, see Lemma
\ref{lemm_stab_psi}.

In addition, for $h^{-1}xh\in P$, we have $\pi^{U_1\rtimes
V'}_{U_1}(h^{-1}xh)=\pi^{P}_{U_1}(h^{-1}xh)$. Using step 2,
we obtain
$$
\begin{array}{ll}
&((\theta\circ f\circ\ln\circ\pi^{U_1\rtimes
V'}_{U_1})\cdot\Ind{P}{U_1\rtimes V'}{(\theta\circ
f\circ\ln\circ\pi^P_{\wt P})})(x)=\\
&=\sum_{h\in H\mid h^{-1}xh\in
P}\theta(f(\ln\pi^{P}_{U_1}(x)+\ln\pi^P_{\wt P}(h^{-1}xh)))=\\
&=\sum_{h\in H\mid h^{-1}xh\in P}\theta(f(\ln
h^{-1}xh))=\Ind{P}{U_1\rtimes V'}{(\theta\circ f\circ\ln)}.
\end{array}
$$

7. Combining steps 5 and 6, we conclude that $$\eta=\Ind{U_1\rtimes
V'}{U}{\Ind{P}{U_1\rtimes V'}{(\theta\circ
f\circ\ln)}}=\Ind{P}{U}{(\theta\circ f\circ\ln)}=\chi,$$ which completes the proof.}

\sect{Characters of maximal dimension}
\label{sect:char_max_dim}
In this section, we give an alternative proof of the formula for the characters of maximal dimension given in \cite{AndreNeto09}. We also give an independent classification of orbits of maximal dimension previously obtained in \cite{AndreNeto06}. Firstly, we need to recall some notations from \cite{IgnatevPetukhov25}.

Let $X$ be a subset of $\Phi^+$. Set $\ut_X:=\bigoplus_{\alpha\in X}\Fp e_\alpha$ with the (Lie) bracket given by the formula $$[e_\alpha, e_\beta]=\begin{cases}[e_\alpha, e_\beta],& {\rm~if~}\alpha+\beta\in X,\\0,&\text{if }\alpha+\beta\notin X\end{cases}$$ for all $\alpha, \beta\in X$
(note that if $\alpha, \beta, \alpha+\beta\in \Phi$, then $[e_\alpha, e_\beta]$ is proportional to $e_{\alpha+\beta})$. Such a bilinear bracket $[\cdot, \cdot ]$ do not always define a Lie algebra but it does define a Lie algebra under the assumption that $X$ is a $C$-pattern (a shorthand for ``combinatorial counterpart of a pattern subgroup'') or $C$-quattern (a shorthand for ``combinatorial counterpart of a quattern subgroup''), see the next definitions. 

\defi{We will say that $X$ is a {\it C-pattern}  if $X$ satisfies the following condition: $$\text{if }\alpha, \beta\in X\text{ and }\alpha+\beta\in\Phi^+\text{ then }\alpha+\beta\in X.$$
This condition is clearly equivalent to the statement that $\ut_X$ is a sub-Lie algebra of $\ut$.}
\defi{ We will say that $X$ is a {\it C-quattern} if $X=X_+\setminus X_-$ for two C-patterns $X_+, X_-$ such that $X_-\subset X_+$ and
$$\text{if }\alpha_-\in X_-,\alpha_+\in X_+~\text{ and }\alpha_-+\alpha_+\in X_+ \text{ then }\alpha_-+\alpha_+\in X_-.$$
This is clearly equivalent to the statement that $\ut_{X_-}$ is an ideal of the Lie algebra $\ut_{X_+}$ and in this case we have $\ut_X\cong \ut_{X_+}/\ut_{X_-}$ where the right hand side is a quotient of Lie algebra by its ideal and thus the left hand side is a Lie algebra.}
From now on we assume that $X$ is always a C-quattern. We denote by $U_X$ the unipotent group with Lie algebra $\ut_X$. Every C-pattern $X$ is a C-quattern for $X_+=X, X_-={\varnothing}$ and for every C-quattern $X$ the data $(\ut_X, [\cdot, \cdot])$ provides a structure of Lie algebra on $\ut_X$. 

Set 
$$Z(X):=\{\alpha\in X\mid (\alpha+X)\cap X={\varnothing}\}.$$ 
Then it is easy to verify that $\nt_{Z(X)}$ coincides with the center of $\nt_X$. In particular, this implies that $Z(X)\ne\varnothing$ for any C-quattern $X$.

\defi{Pick $f\in \ut_X^*$ and consider $Z\subset Z(X)$. We will say that $f$ is {\it $Z$-saturated} if $f(e_\alpha)\ne0$ for all $\alpha\in Z$. We will say that $f$ is {\it saturated} if $f$ is $Z(X)$-saturated. 
We denote the variety of $Z$-saturated linear forms by $\nt^*_{X; Z}$.}
It is clear that each $\ut_{X; Z}^*$ is $U_X$-stable and thus is a collection of coadjoint orbits. We map $\ut^*_{X; Z}$ to $\ut^*$ by the formula $e_\alpha^*\mapsto e_\alpha^*$ and denote by $\underline{\nt}^*_{X; Z}$ the image of this map. 
We will frequently identify $\ut_{X; Z}^*$ with $\underline{\ut}_{X; Z}^*$. 

\defi{Let $Z\subset  Z(X)$ be a subset and let $Y\subset \ut_{X; Z}^*$ be a subvariety. 
We say that $Y$ is a {\em set-section} for $\ut_{X; Z}^*$ if $Y$ intersects each $U_X$-orbit of $\ut_{X; Z}^*$ in exactly one point.}
It is clear that a set-section is a section of the native quotient map $\ut_{X; Z}^*\to \ut_{X; Z}^*/U_X$  but only in set-theoretic terms. 
Also note that all the set-sections which we will consider explicitly will be unions of the form $V(S_1)\sqcup V(S_2)\sqcup\ldots$ for some subsets $S_1, S_2, \ldots\subset\Phi^+$. 

Next, we define the subset $D_{\text{reg}}$ of positive roots as follows:
$$D_{\text{reg}}=\{(2\epsi_1,2\epsi_2,\ldots,2\epsi_n,\epsi_1+\epsi_2,\epsi_2+\epsi_3,\ldots,\epsi_{n-1}+\epsi_n)\}\in\Phi^+$$
(i.e., the subset $D_{\text{reg}}$ consists of the roots corresponding to the squares  $(-1,1),\ldots,(-n,n)$ and $(-2,1),\ldots,(-n,n-1)$, as in Example \ref{ex:squares}). As in Section \ref{main_definitions}, we will consider a set $\xi$ of constants from $\Fp_q^{\times}$; however in this section, it will be defined with a minor correction. More precisely, put $\xi=(\xi_{\beta})_{\beta\in D}$ for an arbitrary set $D$ of positive roots, but allow $\xi_{2\epsi_n}$ to be zero. Now we are ready to formulate and prove a following result.

\theop{Let\label{theo_orb_max_dim} $D\in\Phi^+$ be a maximal orthogonal subset of $D_{\text{reg}}$\textup{,} and $\xi$ be as above\textup{,} then put 
\begin{equation}\label{lin_form_max_dim}
f=f_{D, \xi} = \sum_{\beta\in D}\xi_{\beta}e_{\beta}^*\in\ut^*.
\end{equation}
Then the orbit of $f$ has maximal dimension. Conversely\textup{,} each orbit of maximal dimension contains exact one such a linear form.}{The fact that all such orbits have maximal dimension follows directly from the construction of polarizations (note that maximal possible dimension of a coadjoint orbit in $\ut^*(C_n)$ is known: it equals to $n(n-1)$, see, for example, \cite{AndreNeto06}). 
More precisely, $\dim\Omega_{f_{D,\xi}}=2\codim\pt=|M_D|$, where $\pt$ is a polarization for $f=f_{D,\xi}$, and it is easy to see that $|M_D|$ is always $(n-1)(n-2)/2$ in that case. So, $\dim\Omega_f=n(n-1)$ is maximal. 

To prove the converse, we firstly want to proof that, given $f=(\xi_{\beta})_{\beta\in\Phi^+}\in\ut^*$ and $\delta\prec'\epsi_1+\epsi_2$, where $\delta$ is the biggest root with nonzero constant $\xi_{\delta}$ with respect to the order $\prec'$, the orbit of $f$ has non-maximal dimension. Let us assume that $\delta=\epsi_1+\epsi_3$, i.e., $\xi_{2\epsi_1}=\xi_{\epsi_1+\epsi_2}=0$. Consider the quattern $Y=\Phi^+\backslash\{2\epsi_1,\epsi_1+\epsi_2\}$ and denote by $\wt f$ the projection of $f$ to $\ut_{Y;Z}$ for $Z=\langle e_{\delta}\rangle\subset Z(Y)$. The dimension of the orbit of the form $f$ can be calculated as $\text{rank}(\beta_{f})$ (recall that $\beta_f$ is a quadratic form on $\ut$ defined by $\beta_f(x,y)=f([x,y])$). We want to compare the ranks of the quadratic forms $\beta_f$ and $\beta_{\wt f}$ defined on $\ut$ and $\ut_{Y;Z}$ respectively. Let us write down the matrix $A$ of the form $\beta_f$ (all roots are sorted in the order $\prec'$):
\begin{equation}A=\begin{pmatrix}
    * & \cdots & * & \cdots & * & 0 & 0\\
    \vdots & \ddots & & & \vdots & \vdots & \vdots\\
    * & & \ddots & & * & -cf(2\epsi_1) & 0\\
    \vdots & & & \ddots & \vdots & \vdots & \vdots\\
    * & \cdots & * & \cdots & * & 0 & 0\\
    0 & \cdots & cf(2\epsi_1)& \cdots & 0 & 0 & 0\\
    0 & \cdots & 0 & \cdots & 0 & 0 & 0
\end{pmatrix}.\end{equation}
Here the entries corresponding to roots $\alpha$ and $\delta$ (and so filled by $\beta_f(e_{\alpha},e_{\delta})=\beta_{\wt f}(e_{\alpha},e_{\delta})$) are marked by $*$ if $\alpha+\delta\notin\{2\epsi_1,\epsi_1+\epsi_2\}$. Clearly, $f(2\epsi_1)=0$, so $\text{rank}(\beta_f)=\text{rank}(\beta_{\wt f})$. Hence, we can consider the form $\wt f$ and its orbit in $\ut_{Y;Z}$ instead of $f$ and its orbit in $\ut$. 

Let $Y'$ be defined as the set $Y\backslash\{\delta_1,\cdots,\delta_{2n-4},\beta_1\cdots,\beta_{2n-4}\}$, where $\delta_i\in S^+(\gamma)$ and $\beta_i\in S^-(\gamma)$. Consider the map $\vfi:\ut^*_{Y}\to\ut^*=\ut^*_X$ defined as follows: for $e_{\tau}^*\in\ut_Y^*$,
$$\vfi(e^*_{\tau})=e^*_{\tau}\in\ut^*.$$ According to \cite[Proposition 5.11]{IgnatevPetukhov25},
the set $Y'$ is a quattern and there is a form $f'\in\ut_{Y';Z}$ such that $\vfi(f')$ lying on the orbit of the form $\wt f$ in $\ut_{Y;Z}$ and $$\dim\Omega_{\wt f}=\dim\Omega_{f'}+|\{\delta_1,\ldots,\delta_{2n-4},\beta_1,\ldots,\beta_{2n-4}\}|=\dim\Omega_{f'}+2(2n-4).$$
One can easily see that the quotient algebra $\ut_{Y'}$ is isomorphic to the sum of its ideals $\ut(C_{n-2})\oplus \langle e_{\epsi_1+\epsi_3}\rangle\oplus\langle\epsi_2-\epsi_3\rangle$. It follows that $\dim\Omega_{f'}$ in $\ut_{Y;Z}$ is less than $(n-2)(n-3)$, the maximal possible dimension of orbit in $\ut(C_{n-2})$. Hence, we obtain that $$\dim\Omega_{\wt f}\leq (n-2)(n-3)+2(2n-4)=n(n-1)-2<n(n-1).$$

The other cases when $\delta\prec'\epsi_1+\epsi_3$ can be considered similarly. For example, let $\delta=\epsi_1+\epsi_4$. In that case, the orbit of the form $f$ is contained in the quotient algebra of the algebra $\ut_{Y}$ by the ideal $\langle e_{\epsi_1+\epsi_3}\rangle\subset Z(Y)$, but we have already proved that the maximal possible dimension of an orbit in $\ut_{Y;Z}$ is less than $n(n-1)$, where $Z=\{\epsi_1+\epsi_3\}$. It is left to note that $\ut_{Y;Z}$ is an open subspace of $\ut_Y$ as the set of orbits of maximal dimension, so its intersection with the set of orbits of maximal dimension is non-empty, so the maximal dimension of an orbit in $\ut_Y$ is also less than $n(n-1)$.

Now, we know that a linear form $f=(\xi_{\beta})_{\beta\in\Phi^+}$ having the orbit of maximal dimension has to have nonzero constant $\xi_\delta$, where $\delta$ is the biggest root such that $\xi_{\delta}\ne 0$ with $\delta=2\epsi_1$ or $\epsi_1+\epsi_2$. We aim to present a set-section of the set of orbits of maximal dimension. Let us consider the case when $\delta=2\epsi_1$. We proceed by induction on $n$. The base of induction for $\Phi=C_2$ is known, see, f.e., \cite{Goodwin}. 
Fix $f$ and put $X=\Phi^+$, $\gamma=\delta=2\epsi_1\in Z(X)$, and $Y=X\backslash\{\delta_1,\dots,\delta_{n-1},\beta_1,\dots,\beta_{n-1}\}$ for $\delta_i=\epsi_1-\epsi_i$ and $\beta_i=\epsi_1+\epsi_i$. The form $f$ belongs to $\ut_{X;Z(X)}$ by definition.

The Proposition 5.11 tells us that $Y$ is a quattern and if $K$ is a set-section for $\ut_{Y;Z(X)}$ then $\vfi(K)$ is a set-section for $\ut_{X;Z(X)}$. Note also that $\ut_{Y}$ is isomorphic to $\ut(C_{n-1})\oplus I$, where $I=\langle e_{2\epsi_1}\rangle $ and both of the summand are ideals. It means that we can pick a form $f'\in\ut_{Y,Z(X)}$ such that $\wt f=\vfi(f')$ lies on the orbit of $f$ (note that $\wt f$ has zeroes at places corresponding to the roots from $C_1$). According to the Proposition 5.11 again we also claim that $$\dim\Omega_{f}=\dim\Omega_{\wt f}=\dim\Omega_{f'}+2(n-1)\leq (n-1)(n-2)+2(n-1)=n(n-1),$$ and we have the equality if and only if $\dim\Omega_{f'}=(n-1)(n-2)$, is maximal. Thus, by the inductive assumption, $f'$ has a linear form $g=g_{\wt D,\wt\xi}$ for $\wt D=D'\cup 2\epsi_1$, where $D'$ is a maximal orthogonal subset of $D'_{\text{reg}}\subset\Phi^+_{n-1}$ and $\wt\xi=(\xi_{\beta})_{\beta\in\wt D}$, on its orbit. Hence, $f$ has the element $f_{\wt D,\wt\xi}$ of the form \eqref{lin_form_max_dim} on its orbit. That finishes the proof in this case.

The case when $\delta=\epsi_1+\epsi_2$ can be considered similarly. The main difference is that $X$ must be chosen as $\Phi^+\backslash2\epsi_1$, $\gamma=\delta=\epsi_1+\epsi_2\in Z(X)$ and for $Y$ we take the set $X\backslash\{\delta_1,\ldots,\delta_{2n-3},\beta_1,\ldots,\beta_{2n-3}\}$, where
$\delta_i\in S^+(\gamma)$ and $\beta_i\in S^-(\gamma)$. The next arguments are the same as in the previous case.}

As a corollary, we can compute the number of such orbits.
\corop{Denote by $S_n$ the number of the coadjoint orbits of maximal dimension over the field $\Fp_q$\textup, then 
$S_n=Ar_1^n+Br_2^n$\textup, where
\begin{equation*}
    A=\frac{v+2+\sqrt{v^2+4v}}{2\sqrt{v^2+4v}},\quad B=\frac{-v-2+\sqrt{v^2+4v}}{2\sqrt{v^2+4v}}
\end{equation*}
and
\begin{equation*}
    r_1=\frac{v+\sqrt{v^2+4v}}{2},\quad r_2=\frac{v\sqrt{v^2+4v}}{2}
\end{equation*} for $v=q-1$. In particular\textup, $S_n$ is a polynomial with integer non-negative coefficients.}{From the proof of Theorem \ref{theo_orb_max_dim} we see that the sequence $S_n$ has a simple recursive condition $S_n=vS_{n-1}+vS_{n-2}$. Solving this linear recurrent, we obtain the requested result.}

\exam{Here we schematically drew an example of a form of maximal dimension. The set $D$ equals $\{\epsi_1+\epsi_2,2\epsi_3,2\epsi_4\}$, the symbol $\otimes$ denotes a nonzero constant, and the symbol $\square$ denotes an arbitrary constant from the field $\Fp_q$.
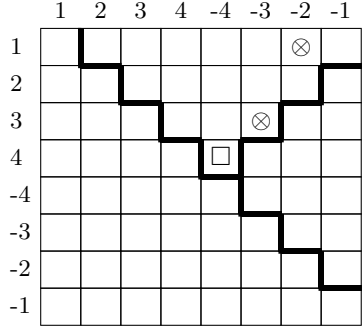
\begin{center}\small
$\mymatrix{
\lNote{1}\Note{1}\pho& \Lft{2pt}\Bot{2pt}\Note{2}\pho& \Note{3}\pho& \Note{4}\pho& \Note{-4}\pho& \Note{-3}\pho& \Note{-2}\otimes& \Bot{2pt}\Note{-1}\pho\\
\lNote{2}\pho& \pho& \Lft{2pt}\Bot{2pt}& \pho& \pho& \pho& \Bot{2pt}\Rt{2pt}& \pho\\
\lNote{3}\pho& \pho& \pho& \Lft{2pt}\Bot{2pt}& \pho& \Bot{2pt}\Rt{2pt}\otimes& \pho& \pho\\
\lNote{4}\pho& \pho& \pho& \pho& \Rt{2pt}\Lft{2pt}\Bot{2pt}\square& \pho& \pho& \pho\\
\lNote{-4}\pho& \pho& \pho& \pho& \pho& \Lft{2pt}\Bot{2pt}& \pho& \pho\\
\lNote{-3}\pho& \pho& \pho& \pho& \pho& \pho& \Lft{2pt}\Bot{2pt}\pho& \pho\\
\lNote{-2}\pho& \pho& \pho& \pho& \pho& \pho& \pho& \Lft{2pt}\Bot{2pt}\\
\lNote{-1}\pho& \pho& \pho& \pho& \pho& \pho& \pho& \pho\\
}$
\end{center}
}

Next, we need to define another important subsets of the positive roots. Let $D\subset D_{\text{reg}}$ be a maximal orthogonal subset. Consider the set $I_D$ of the natural numbers $I=\{\col(\beta)\mid\beta=\epsi_i+\epsi_j,\beta\in D\}$ and put $D'=D\bigcup_{i\in I_D}2\epsi_i$.

\defi{Let $D\subset\Phi^+$ be an arbitrary subset of positive roots and $\phi=(\phi_{\beta})_{\beta\in D}$ be a set of nonzero constant from the field $\Fp_q$. Define the element $x_D(\phi)\in\ut$ as follows:
$$x_D(\phi)=\sum_{\beta\in D}\phi_{\beta}e_{\beta}.$$}

\defi{Let $D\in\Phi^+$ and $\beta\in\Phi^+$. Denote 
$$X_{\beta}^D=\{\row(\alpha)\mid\alpha\in D, \row(\alpha)\prec\row(\beta),\col(\alpha)\succ\col(\beta)\}\cup\row(\beta),$$
$$Y_{\beta}^D=\{\col(\alpha)\mid\alpha\in D, \row(\alpha)\prec\row(\beta),\col(\alpha)\succ\col(\beta)\}\cup\col(\beta)$$
sorted by the order $\prec$. Denote the minor with the rows from $X_{\beta}^D$ and the columns from $Y_{\beta}^D$ by $\Delta_{\beta}^D$.
}

The next step before we can formulate the main result of this section is to define certain varieties. For given $D\subset\Phi^+$ and $\phi
\colon D\to\Fp_q^{\times}$ we define the variety
$$k_D(\phi)=\{x\in\ut\mid\Delta_{\beta}^D(x)=\Delta_{\beta}^D(x_D(\phi))\text{ for all }\beta\in R(D)\}.$$
Here $R(D)=\Phi^+\backslash S(D)$ for $S(D)=\bigcup_{\alpha\in D}(\alpha+\Phi^+)$ (by definition, the set $\alpha+\Phi^+$ contains all the roots $\alpha+\gamma,\gamma\in\Phi^+$, where this sum is also a positive root). We will also write $K_D(\phi)=\exp(k_D(\phi))$ and $X_D(\phi)=\exp(x_D(\phi))$.

Now, we recall the linear character $\eta$ of the group $\Fp_q^{\times}$ defined by
$$\eta(c)=\begin{cases}
    1,\text{ if }c\in(\Fp_q^{\times})^2,\\
    -1, \text{ otherwise},
\end{cases}$$
for all $c\in\Fp_q$. Furthermore, for the characters $\eta$ and $\theta$, we define the Gauss sum of $\eta$ and $\theta$ as 
$$G(\eta,\theta)=\sum\limits_{c\in\Fp_q^{\times}}\eta(c)\theta(c).$$

For the convenience of writing, we introduce some additional notations. Let $j$ be 0 or 1, put
\begin{equation*}
D_j=
\begin{cases}
D\cap \{\epsi_i+\epsi_{i+1}\mid1\leq i\leq n-1\},&\text{ if }j=1,\\
D\cap\{2\epsi_i\mid1\leq i\leq n\},&\text{ if }j=0,
\end{cases}
\end{equation*}
and
\begin{equation*}
D_j(\alpha)=
\begin{cases}
D\cap\{\epsi_i+\epsi_{i+1}\mid \epsi_i+\epsi_{i+1}\prec'\alpha\},&\text{ if }j=1,\\
D\cap\{2\epsi_i\mid 2\epsi_i\prec'\alpha\},&\text{ if }j=0,\\
\end{cases}
\end{equation*}
whenever $\alpha\in\Phi^+$. 

%Define also a set $D(\alpha)=\bigcup_{\beta\in T(D,\alpha)}\Et(\beta)$, where $T(D,\alpha)=\{\gamma\in D\mid\row(\gamma)<\row(\alpha),\col(\gamma)>\col(\alpha)\}$ for $\alpha\in\Phi^+$. And \begin{equation*}\Et(\beta)=
%\begin{cases}
%(j,i)\cup (-i,-j)\text{ if }\beta=\epsi_i-\epsi_j,\\
%(-j,i)\cup (-i,j)\text{ if }\beta=\epsi_i+\epsi_j,\\
%(-i,i)\text{ if }\beta=2\epsi_i.
%\end{cases}
%\end{equation*}

Now we are ready to describe the support of the character corresponding to an orbit of maximal dimension and calculate its value on a conjugacy class from the support.

\mtheo{Let\label{theo_char_value} $D$ be a maximal orthogonal subset of $D_{reg}$\textup, $\xi$ as above\textup, and $f=f_{D,\xi}$. Denote by $\chi$ the character corresponding to the orbit $\Omega_{D,\xi}$ of the form $f$\textup, then
\begin{enumerate}
    \item [\textup{a) }]$\Supp{\chi}=\bigsqcup\limits_{D^*,\phi}K_{D^*}(\phi)$\textup, where $D^*$ runs over all possible orthogonal subsets of $D'$ and $\phi$ runs over all maps from $D^*$ to $\Fp_q^{\times}$.
    \item [\textup{b) }]$\chi(K_{D^*}(\phi))=C\prod\limits_{\alpha\in D^*\cap D}\theta(\xi_{\alpha}\phi_{\alpha}),$
\end{enumerate}
where $$C=q^{(\frac{n(n-1)}{2}-t(D,D^*)-t_1(D,D^*))} G(\eta,\theta)^{t_0(D,D^*)}\prod\limits_{\alpha\in D_0}\prod\limits_{\beta\in D^*_0(\alpha)}\eta(\xi_{\alpha}\phi_{\beta}),$$ and\textup, by definition\textup,
\begin{equation*} t_0(D,D^*)=\sum\limits_{\alpha\in D_0}|D_0(\alpha)|,\quad t_1(D,D^*)=\sum\limits_{\alpha\in D_1}|D_1^*(\alpha)|,\quad
t(D,D^*)=\sum\limits_{\alpha\in D}|D_0^*(\alpha)|+|D_1^*(\alpha)|.
\end{equation*}}
 
In fact, in the process of proving this theorem, we will obtain that $$\ln(\Supp{\chi})=\bigsqcup_{D^*,\phi}k_{D^*}(\phi),$$ and the varieties $k_{D*}(\phi)$ coincide with the orbits of the adjoint action of elements $x_{D^*}(\phi)$, but $\exp(gxg^{-1})=g\exp(x)g^{-1}$ for any $x\in\ut$. So, the variety $$K_{D^*}(\phi)=\exp(k_{D^*}(\phi))=\{\exp(gx_{D^*}(\phi)g^{-1})\mid g\in U\}$$ coincides with the conjugacy class of $X_D(\phi)$. Among other things, it gives us that the notation $\chi(K_{D^*}(\phi))$ is correct. Before proving Theorem \ref{theo_char_value}, we need to prove a technical lemma.

\lemmp{Let $D$ be as in Theorem \textup{\ref{theo_char_value}} and $D^*\in D'$ be an arbitrary rook placement, put $y=x_{D^*}(\phi)$ where $\phi:D^*\to\Fp_q^{\times}$. Then $|U.y|=|k_{D^*}(\phi)|$\textup, where $U.y$ denotes the orbit of $y$ under the adjoint action of the group $U$.
}{Firstly, denote the variety $k_{D^*}(\phi)$ by $Y$, for convenience, and notice that $Y$ is isomorphic to the affine space over the filed $\Fp_q$ of dimension $S(D)$ (so its cardinality is equal to $q^{S(D)}$). The short explanation of this statement is that equations of the form $\Delta_{\beta}^{D^*}(x)=\Delta_{\beta}^{D^*}(y)$ on the elements $x$ from $Y$ are in fact the equations of the form $x_{i,j}=f_{i,j}(x_{-1,1},x_{-2,1},\ldots,\widehat{x_{i,j}},\cdots)$, what gives us a regular bijective map between $Y$ and the affine space of the requested dimension, cf. \cite{Ignatev09}.
Using the fact that $|U.x|=\dfrac{\mid U\mid}{|\Stab~y|}$, we can calculate the cardinality of the orbit of the element $y$ by calculating its stabilizer. Indeed, for $g=\exp z\in U$, where $z=\sum_{\alpha\in\Phi^+}z_{\alpha}e_{\alpha}$, we have
\begin{equation}\label{form_ad_action}g.y=gyg^{-1}=(\exp\ad{z})(y)=y+[z,y]+\frac{1}{2}[z,[z,y]]+\ldots.\end{equation}
By $C(\beta)$ we will denote the set of positive roots whose sum with $\beta$ is also a positive root. It is easy to see that
$$C(\beta)=\begin{cases}
\{\epsi_j-\epsi_i\mid j<i\}\text{ if }\beta=2\epsi_i,\\
\{\epsi_j-\epsi_i\mid j<i\}\cup\{\epsi_{i+1}-\epsi_j\mid j<i+1\}\text{ if }\beta=\epsi_i+\epsi_{i+1}.
\end{cases}$$
Notice that $C(\beta_1)\cap C(\beta_2)=\varnothing$ for any $\beta_1$ and $\beta_2$ from $D^*$. Next, denote by $\delta$ the smallest root from $D^*$ with respect to the order $\prec'$, and denote the set $\cup_{\beta\in D^*}C(\beta)$ by $C(D^*)$. Now, denote by $\delta'$ the smallest root with respect to the same order from all $\alpha\in C(\delta)$ such that $z_{\alpha}\ne 0$. From formula \eqref{form_ad_action} and the fact that $\beta\prec'\alpha+\beta$ for $\beta\in D^*$ and $\alpha\in\Phi^+$ we see that $(g.y)_{\delta+\delta'}=z_{\delta'}$, because the root $\delta+\delta'$ is the closest root to $\delta$ which is smaller than $\delta$. But $y_{\delta+\delta'}=0$, so we conclude that $z_{\delta'}$ must be zero what contradicts to the existence of $\delta'$. So, $z_{\alpha}=0$ for all $\alpha\in C(\delta)$. Now we can apply the same arguments to the next smallest root from $D^*$. Hence, after some iterations we get that $z_{\alpha}=0$ for all $\alpha\in C(D^*)$. The last thing we need to notice is that $|C(D)|=|S(D)|$. Thus, $\ln(\Stab~y)$ is isomorphic to the affine space of dimension $\dim\ut-|S(D)|$, so $|U.y|=q^{|S(D)|}$. This concludes the proof.
}

Notice also that $U.y$ is a subset of $k_{D^*}(\phi)$ for $y=x_{D^*}(\phi)$ for $D^*$ as in the lemma above. It follows from the statement that $\Delta^{D^*}_{\beta}(gyg^{-1})=\Delta^{D^*}_{\beta}(y)$. Thus, the sets $U.y$ and $k_{D^*}(\phi)$ coincide.

\medskip\textsc{Proof of Theorem~\ref{theo_char_value}, part ({\normalfont a})}.

We will start from the first statement of the theorem. We already have that the varieties $K_{D^*}(\phi)$ coincide with the conjugacy classes of the elements $X_{D^*}(\phi)$. Firstly, we will split the support of $\chi$ into the bigger varieties $\wt K_{D^*}(\phi)$. Recall the definitions of $\chi$, $U_1$, $V'$ and $\wt U$ from the previous section. We proceed by induction on $n$. The base of induction can be checked trivially.

Assume that $\beta_1=2\epsi_1$ and put $g_{U_1}=\pi_{U_1}^{U_1\rtimes V'}$ and $g_{\wt U}=\pi_{\wt U}^{U_1\rtimes V'}$ for an arbitrary $g=\exp(x)\in U$. Denote also a complete set of coset representatives of $U_1\rtimes V'$ (will write it just as $U_1V'$, for convenience) in $U$ by $M$. Using the standard formula for the induced character, we obtain that
\begin{equation}\label{form_char}
\chi(g)=\sum\limits_{y=h^{-1}gh\mid y\in U_1V',h\in M}\theta(f(\ln(y_{U_1})))\wt\chi(y_{\wt U})
\end{equation}
(recall that in that case $\wt U=V'$ and $\wt\ut$ is isomorphic to $\ut(C_{n-1})$). More precisely, $M$ can be chosen as follows: $$M=\{\exp(t)\mid t\in \ut\backslash(\ut_1\oplus\vt)\}.$$

The condition $y\in U_1V'$ might be reformulated as $h^{-1}xh\in\ut_1\oplus\vt$. Indeed, $y=h^{-1}\exp(x)h=\exp(h^{-1}xh)$, but $U_1V'=\exp\ut_1\exp\vt\subset\exp(\ut_1\oplus\vt)$ obviously, and these sets have the same cardinality, so they coincide. 

The conjugation by elements $h\in M$ does not affect the elements that are not in $C_1$ and $R_{-1}$ in~$g$. The projection $\pi$ on the group $\wt U$ works by deleting the elements in the first column, so $y_{\wt U}=g_{\wt U}$. Hence, formula \eqref{form_char} can be written as
$$\chi(g)=\wt\chi(g_{\wt U})\sum\limits_{y=h^{-1}gh\mid y\in U_1V',h\in M}\theta(f(\ln(y_{U_1}))).$$
Taking into account that $g_{\wt U}$ must belong to $\Supp{\wt\chi}$, we obtain a part of the equations on $g$ that we need by applying the induction assumption to it. More precisely,
$$\Delta_{\beta}^{D^*}(x)=\Delta_{\beta}^{D^*}(x_{D^*}(\vfi))\text{ for all }\beta\in R(D^*)\backslash C_1$$
for certain $D^*\in D'$. Combining these conditions and the condition that $h^{-1}xh\in\ut_1\oplus\vt$ we also obtain that
$$\Delta_{\beta}^{D^*}(x)=\Delta_{\beta}^{D^*}(x_{D^*}(\vfi))\text{ for all }\beta\in C_1\backslash C'.$$
In fact, all these equations in that case mean only that $x_{1,i}=0$ for all $i$ from 1 to $n$. Thus, it remains to obtain the required equations for the roots from the set $C'$.

To get that equations, we need to examine the expression $\theta(f(\ln(y_{U_1})))$. Because in that case $f(x)=x_{-1,1}$ for $x\in\ut_1$ we only interested in the element $(\ln(y_{U_1}))_{-1,1}$. Straightforward matrix calculation shows us that
$$\theta(f(\ln(y_{U_1})))=\theta(\xi_{\beta_1}(\sum\limits_{i,j=2}^{n}h_{i,1}h_{j,1}x_{-i,j}-2\sum\limits_{i=2}^nh_{i,1}x_{-i,1}+x_{-n,1}))$$
(we remind that $x_{-i,j}=x_{-j,i}$, where $i=1,\ldots,n$). We see that the expression in the parentheses is a quadratic form over the finite field in the variables $h_{i,1}$, so we can reduce it to the canonical form. 

There are two possibilities:
\begin{enumerate}
    \item This form is not central; therefore, its canonical form will be
    $\lambda_1z_1^2+\lambda_2z_2^2+\ldots+\lambda_kz_k^2+\lambda_{k+1}z_{k+1}.$
    \item This form is central, so its canonical form will be
    $\lambda_1z_1^2+\lambda_2z_2^2+\ldots+\lambda_kz_k^2.$
\end{enumerate}
(where $\lambda_i$ are some nonzero constants and $z_j$ are variables). According to the fact that $\sum_{c\in\Fp_q}\theta(c)=0$ we conclude that the first option is bad for us because the sum over all $h$ will be zero (due to the linear summand $\lambda_{k+1}z_{k+1}$), so the value $\chi(g)$ will be zero, too. Thus, our form has to be central. We can write down the matrix of this quadratic form. It is easy to check that it will be exactly the matrix consisting of the elements $\wt x_{i,j}$ of the matrix $\pi(x)$ (which are, in fact, elements $x_{i,j}$ of the matrix $x$ without the first column and the last row); let us denote this matrix by $\wt x$. The coefficients of the linear part of this form are $(-2x_{-2,1},-2x_{-3,1},\ldots,-2x_{-n,1}$); let us denote this vector by $b$. The set of all centres of this quadratic form is described by the system of linear equations $\wt xy^t=-\frac{b^t}{2}$, where $y$ is the vector of variables. Applying the Kronecker--Capelli criterion to this system, we obtain the remaining part of equations on the first column, except the conditions on the minor $\Delta_{2\epsi_1}^{D^*}(x)$, what gives us the variety $\wt k_{D^*}(\phi)$ where this minor might have any value. It is equal to the condition that the root $\beta_1=2\epsi_1$ can belong to $D^*$ or not. Hence, $x\in k_{D^*}(\phi)$ and $g\in K_{D^*}(\phi)$ for some $D^*\in D'$. 

Let us move to the case $\beta_1=\epsi_1+\epsi_{n-1}$. As in the previous case, we will start from describing a complete set of cosets representatives of $U_1V'$. One can easily check that set $M$ can be chosen in the following way:
$$M=\{\exp(t)\mid t\in\ut\backslash(\ut_1\oplus\vt')\}.$$

Next, we need to notice that elements not from $C_1$, $C_2$ and $R_{-1}$, $R_{-2}$ in arbitrary $g\in U$ do not change by the conjugation with $h\in M$. Furthermore, the projections from $U_1V'$ to $V'$ and from $V'$ to $\wt U$ work just by deleting the elements from $C_1$, $R_{-1}$ and $C_{2}$, $R_{-2}$ respectively. So the formula \eqref{form_char} can be rewritten as 
$$\chi(g)=\wt\chi(g_{\wt U})\sum\limits_{y=h^{-1}gh\mid y\in U_1V',h\in M}\theta(f(\ln(y_{U_1}))).$$

Since $\wt\chi(g_{\wt U})\ne0$, we obtain the following equations on the element $x$:
$$\Delta_{\beta}^{D^*}(x)=\Delta_{\beta}^{D^*}(x_{D^*}(\vfi))\text{ for all }\beta\in R(D^*)\backslash (C_1\cup C_2)$$
for certain $D^*\in D'$. Now, we want to satisfy the condition $y=h^{-1}gh\in U_1V'$. As in the previous case, this statement is equivalent to the statement that $z=h^{-1}xh\in\ut_1\oplus\vt'$ by arguments as above. From the conditions $z_{i,2}=0$ for $j=\overline{3,n}$ and the equations on the minors written above, we obtain that $x_{j,2}=0$ for $j=\overline{3,n}$ (for this places $\Delta_{\epsi_2-\epsi_i}^{D^*}(x)=x_{i,2}$). Let us write down the remaining conditions which appear from the condition that $z\in\ut_1\oplus\vt'$. The conditions $z_{i,2}=0$ for $i$ from $-n$ to $-3$ give us that $\wt x\wt h^t=b^t$, where $\wt x$ is $(2n-2)\times(2n-2)$ matrix consists of the elements of the matrix~$\pi(x)$, $\wt h=(h_{3,2},h_{4,2},\ldots,h_{n,2})$ and $b=(h_{-n,2},h_{-n+1,2},\ldots,h_{-3,2})$. The condition $z_{-2,2}=0$ gives us
$$\sum\limits_{i,j=3}^{n}h_{i,2}h_{j,2}x_{-i,j}-2\sum\limits_{i=3}^nh_{i,2}x_{-i,2}+x_{-n,2}=0.$$
Using equations from the system $\wt x\wt h^t=b^t$ we can rewrite latter expressions as follows:
$$\sum\limits_{i=3}^{n}x_{-i,2}h_{i,2}=x_{-2,2}.$$
Combining this equation with the system $\wt x\wt h^t=b^t$ we have the new system of linear equations $\wt{\wt x}\wt h^t=\wt b^t$, where $\wt{\wt x}$ is obtained from $\wt x$ by adding one string $(x_{-3,2},x_{-4,2},\ldots,x_{-n,2})$ and the vector $\wt b$ is obtained from $b$ by adding $x_{-2,2}$. This system also have to be consistent, so we can use Kronecker--Capelli criterion. It gives us the remaining conditions on the second column, in particular the minor $\Delta_{2\epsi_2}^{D^*}(x)$ is always zero (because the rank of the matrix $\wt x$ is not bigger than $n-2$, and the size of the augmented matrix of the system $\wt{\wt x}\wt h^t=\wt b^t$ is $(n-1)\times (n-1)$). Hence, now we have the following equations:
$$\Delta_{\beta}^{D^*}(x)=\Delta_{\beta}^{D^*}(x_{D^*}(\vfi))\text{ for all }\beta\in R(D^*)\backslash C_1$$

To obtain the equations for the first column we have to look at the expression $$\theta(\xi_{\beta_1}(\ln(y_{U_1})))=\sum\limits_{i=2}^{n}x_{i,1}h_{-i,1}-\sum\limits_{i=2}^{n}x_{-i,1}h_{i,1}+x_{-1,1}.$$
From the system $\wt xh^t=b^t$ we see that some of the elements $h_{i,2}$ (and only them), $i=3,\ldots,n$, can be expressed from the others. As a corollary from the expression above, we immediately obtain that $x_{i,1}=0$ for $i=2,\ldots,n$. Hence,
\begin{equation}\label{eq1}
\theta(\xi_{\beta_1}(\ln(y_{U_1})))=-\sum\limits_{i=2}^{n}x_{-i,1}h_{i,1}+x_{-1,1}.
\end{equation}
So, we need to solve the system of linear equations $\wt xh^t=b^t$. Let us consider two sets of indices sorted with respect to the order $\prec$:
$$T_{D^*}=\{\row(\alpha)\mid\alpha\in D^*,\row(\alpha)<-2\}\cup\{-\col(\alpha)\mid\alpha\in D^*,\col(\alpha)>2\},$$
$$T^{D^*}=\{-\row(\alpha)\mid\alpha\in D^*,\row(\alpha)<-2\}\cup\{\col(\alpha)\mid\alpha\in D^*,\col(\alpha)>2\}.$$
Denote the matrix lying on the intersection of the strings from $T_{D*}$ and the columns from $T^{D^*}$ by $x'$. This matrix will be a matrix of maximal order consisting of the elements from $\wt x$ with nonzero determinant, denote this determinant by $\Delta$. After solving the system $\wt xh^t=b^t$ and substituting the result into the equality \eqref{eq1} we obtain that 
$$\theta(\xi_{\beta_1}(\ln(y_{U_1})))=-\sum\limits_{i\in \{1,\ldots,n\}\backslash T^{D^*}}\Delta_{-i,1}h_{i,1}+c,$$
here the minors $\Delta_{-i,1}$ are lying on the strings from $T_{D^*}\cup \{-i\}$ and the columns $T^{D^*}\cup \{1\}$ and $c$ is a certain constant. Thus, all the minors $\Delta_{-i,1}$ for $i\in \{1,\ldots,n\}\backslash T^{D^*}$ must be zero as the coefficients before linear summands. It gives us the remaining part of the equations on the first column, except for the minors $\Delta_{\epsi_1+\epsi_{n-1}}^{D^*}(x)$ and $\Delta_{2\epsi_1}^{D^*}(x)$, that gives us a variety $\wt k_{D^*}(\phi)$, where these two minors might have any value. It is equivalent to the following: the root $\epsi_1+\epsi_{n-1}$ belongs to $D^*$, or $2\epsi_1$ belongs to $D^*$ and $\epsi_1+\epsi_{n-1}$ does not belong to $D^*$, or both of these roots do not belong to $D^*$. That concludes the proof of the first statement of the theorem. \hfill$\square$

\medskip\textsc{Proof of Theorem~\ref{theo_char_value}, part ({\normalfont b})}.

Given $D$, $D^*\subset D'$, $\phi$ and $\xi$ as in the theorem. We will consider the two cases. At first, we will consider the case when $\beta_1=2\epsi_1\in D$. For convenience, we will write $\phi_{\alpha}=0$ if $\phi$ is not defined on $\alpha\in\Phi^+$. We proceed by induction on $n$. The base of induction in this case for $\Phi=C_2$ can be checked directly by the formula $$\chi=\chi_{\Omega_{f}}(g)=q^{-\frac{1}{2}\dim\Omega_f}\sum\limits_{f'\in\Omega_f}\theta(f'(\ln(g))),$$
where $f=f_{D,\xi}$ and $g=X_{D^*}(\phi)=\exp(x_{D^*}(\phi))$. For example, let us consider this calculation for $\beta_1=2\epsi_1$ (in fact, this is the most technically complicated case). We have
$$\chi(g)=q^{-1}\sum\limits_{f'\in\Omega_f}\theta(f'(x_{D^*}(\phi)).$$
According to the formula of the coadjoint action and noticing that $f'(x_{D^*}(\phi))$ is taken $|\Stab(f)|$ times we obtain that
$$\chi(g)=q^{-1}\theta(\xi_{2\epsi_2}\phi_{2\epsi_2})\theta(\xi_{2\epsi_1}\phi_{2\epsi_1})\sum\limits_{y\in\ut}\theta(\xi_{2\epsi_1}(y_{2,1}^2\phi_{2\epsi_2}))|\Stab(f)|^{-1}.$$
Using \cite[Theorem 5.3.3]{LidlNiederreither}
and taking into account that $|\Stab(f)|=|U|/|\Omega_f|$ we get
\begin{equation*}\chi(g)=
\begin{cases}
q\prod\limits_{\alpha\in D\cap D^*}\theta(\xi_{\alpha}\phi_{\alpha}),&\text{ if }\phi_{2\epsi_2}\notin D^*,\\
\eta(\xi_{2\epsi_1}\phi_{2\epsi_2})G(\eta,\theta)\prod\limits_{\alpha\in D\cap D^*}\theta(\xi_{\alpha}\phi_{\alpha})&\text{ otherwise},
\end{cases}
\end{equation*}
which agrees with the formula. Thus, the base of induction is proven in that case. Other cases can be checked similarly.

So, let us move to the step of induction. Recall the definitions of $M$, $U_1$, $\wt U$, $V'$ from the proof of the previous part. Since our character can be obtained as an induced character, we can use the standard formula for an induced character:
$$\chi(g)=\wt\chi(g_{\wt U})\sum\limits_{y=h^{-1}gh|y\in U_1V',h\in M}\theta(f(y_{U_1})).$$
Here $g_{\wt U}$ is the element $x_{\wt D^*}(\wt\phi)\in U(C_{n-1})$, where $\wt D^*=\pi(D^*)$ and $\wt\phi$ is the natural restriction of $\phi$ to the set $\wt D^*$. Proceeding by the formula of the coadjoint action, we obtain
$$\chi(g)=\wt\chi(g_{\wt U})\theta(\xi_{2\epsi_1}\phi_{2\epsi_1})\sum\limits_{y=h^{-1}gh|y\in U_1V',h\in M}\theta(\xi_{2\epsi_1}(\sum\limits_{\alpha\in\wt D^*_0}\phi_{\alpha}h_{\alpha}^2+\sum\limits_{\beta\in\wt D^*_1}2\phi_{\beta}h_{\beta}h_{\beta'})),$$
where $h_{\alpha}=h_{i,1}$ for $\alpha=2\epsi_i\in\wt D^*_0$ and $h_{\beta}=h_{k,1}$, $h_{\beta'}=h_{l,1}$ for $\beta=\epsi_k+\epsi_l\in\wt D^*_1$. Using the argumentation as above, we see that
$$\chi(g)=\wt\chi(g_{\wt U})\theta(\xi_{2\epsi_1}\phi_{2\epsi_1})G^{|\wt D^*_0|}(\eta,\theta)q^{(n-1)-|\wt D^*|}\prod\limits_{\beta\in\wt D^*_0}\eta(\xi_{2\epsi_1}\phi_{\beta}).$$ Taking into account that $D_0=\wt D_0\cup\{2\epsi_1\}$, $D_1=\wt D_1$ and
\begin{equation*}
    t(D,D^*)=t(\wt D,\wt D^*)+|\wt D^*|,\quad t_0(D,D^*)=t_0(\wt D,\wt D^*)+|\wt D^*_0|,\quad t_1(D,D^*)=t_1(\wt D,\wt D^*),
\end{equation*}
and by applying the inductive assumption to $\wt\chi(g_{\wt U})$, we obtain the requested formula.

In fact, the case when $\beta_1=\epsi_1+\epsi_2$ is much simpler. We will skip the proof for the base of induction because it is similar to the previous case. Here we have
$$\chi(g)=\wt\chi(g_{\wt U})\sum\limits_{y=h^{-1}gh|y\in U_1V',h\in M}\theta(f(y_{U_1})).$$
Simple calculations with matrices show us that
$$\chi(g)=\wt\chi(g_{\wt U})\sum\limits_{y=h^{-1}gh|y\in U_1V',h\in M}\theta(\xi_{\epsi_1+\epsi_2}\phi_{\epsi_1+\epsi_2}).$$
The set of $h\in M$ such that $y\in U_1V'$ can be described as the set of $h\in M$ such that $h_{i,2}=0$ for\break $i\in\{-n,\ldots-2\}\cup T^{D^*}$. It is easy to see that $|T^{D^*}|=| D^*_0(\epsi_1+\epsi_2)|+2| D^*_1(\epsi_1+\epsi_2)|$. Hence,
$$\chi(g)=\wt\chi(g)\theta(\xi_{\epsi_1+\epsi_2}\phi_{\epsi_1+\epsi_2})q^{n-1+n-2-|D^*_0(\epsi_1+\epsi_2)|+2|D^*_1(\epsi_1+\epsi_2)|}.$$
It is only left to notice that $\wt D_0=D_0$ and 
\begin{equation*}
     t(D,D^*)=t(\wt D,\wt D^*)+|\wt D^*|,\quad t_0(D,D^*)=t_0(\wt D,\wt D^*),\quad t_1(D,D^*)=t_1(\wt D,\wt D^*)+|D_1^*(\epsi_1+\epsi_2)|.
\end{equation*}
After applying the inductive assumption to $\wt\chi(g_{\wt U})$ we obtain the required result.\hfill$\square$

\begin{comment}

\end{comment}

%\vspace{5mm}
%\noindent Поступила в редакцию 24/{\it V}/2008;\\
%\noindent в окончательном варианте --- 16/{\it VI}/2008.
%\vspace{5mm}

%\begin{center}
%\maintitle[Communicated by Dr. Sci. (Phys. \& Math.) Prof.
%V.E.\,Voskresenskii] {BASIC SUBSYSTEMS OF ROOT SYSTEMS OF TYPES
%$B_n$ AND $D_n$\\AND ASSOCIATED COADJOINT ORBITS}
%\authorright[2008]
%\anauthor[Ignatev Mikhail Viktorovich
%\email{mihail\_ignatev@mail.ru}, Dept. of Algebra and Geometry,
%Samara State University, Samara, 443011, Russia.]{M.V.\,Ignatev}{}
%\end{center}

\bigskip\textsc{Mikhail Venchakov: National Research University Higher School of Economics, Pokrovsky Boulevard 11, 109028, Moscow, Russia}

\emph{E-mail address}: \texttt{mihail.venchakov@gmail.com}

\end{document}